\documentclass[12pt]{article}
\usepackage{amssymb}
\usepackage{amsmath}
\usepackage{amsthm}
\usepackage{color}
\usepackage{comment}
\usepackage{geometry}		
\usepackage{graphicx}

\let\originalleft\left
\let\originalright\right
\renewcommand{\left}{\mathopen{}\mathclose\bgroup\originalleft}
\renewcommand{\right}{\aftergroup\egroup\originalright}

\begin{document}

\def\b0{{\bf 0}}
\def\bx{{\bf x}}
\def\by{{\bf y}}
\def\bX{{\bf X}}
\def\ee{\varepsilon}
\def\cD{\mathcal{D}}
\def\cO{\mathcal{O}}
\def\cP{\mathcal{P}}
\def\cS{\mathcal{S}}
\def\cX{\mathcal{X}}
\def\cY{\mathcal{Y}}
\def\cZ{\mathcal{Z}}
\def\re{{\rm e}}
\def\ri{{\rm i}}

\newcommand{\removableFootnote}[1]{}

\newtheorem{theorem}{Theorem}[section]
\newtheorem{lemma}[theorem]{Lemma}
\newtheorem{proposition}[theorem]{Proposition}



\title{
Grazing-sliding bifurcations creating infinitely many attractors.
}
\author{
D.J.W.~Simpson\\\\
Institute of Fundamental Sciences\\
Massey University\\
Palmerston North\\
New Zealand
}
\maketitle


\begin{abstract}

As the parameters of a piecewise-smooth system of ODEs are varied,
a periodic orbit undergoes a bifurcation when it collides with a surface where the system is discontinuous.
Under certain conditions this is a grazing-sliding bifurcation.
Near grazing-sliding bifurcations structurally stable dynamics are captured by piecewise-linear continuous maps.
Recently it was shown that maps of this class can have infinitely many asymptotically stable periodic solutions of a simple type.
Here this result is used to show that at a grazing-sliding bifurcation an asymptotically stable periodic orbit
can bifurcate into infinitely many asymptotically stable periodic orbits.
For an abstract ODE system the periodic orbits are continued numerically
revealing subsequent bifurcations at which they are destroyed.

\end{abstract}

\section{Introduction}
\label{sec:intro}
\setcounter{equation}{0}

\begin{figure}[t!]
\begin{center}
\setlength{\unitlength}{1cm}
\begin{picture}(17,4.2)
\put(0,0){\includegraphics[height=4.2cm]{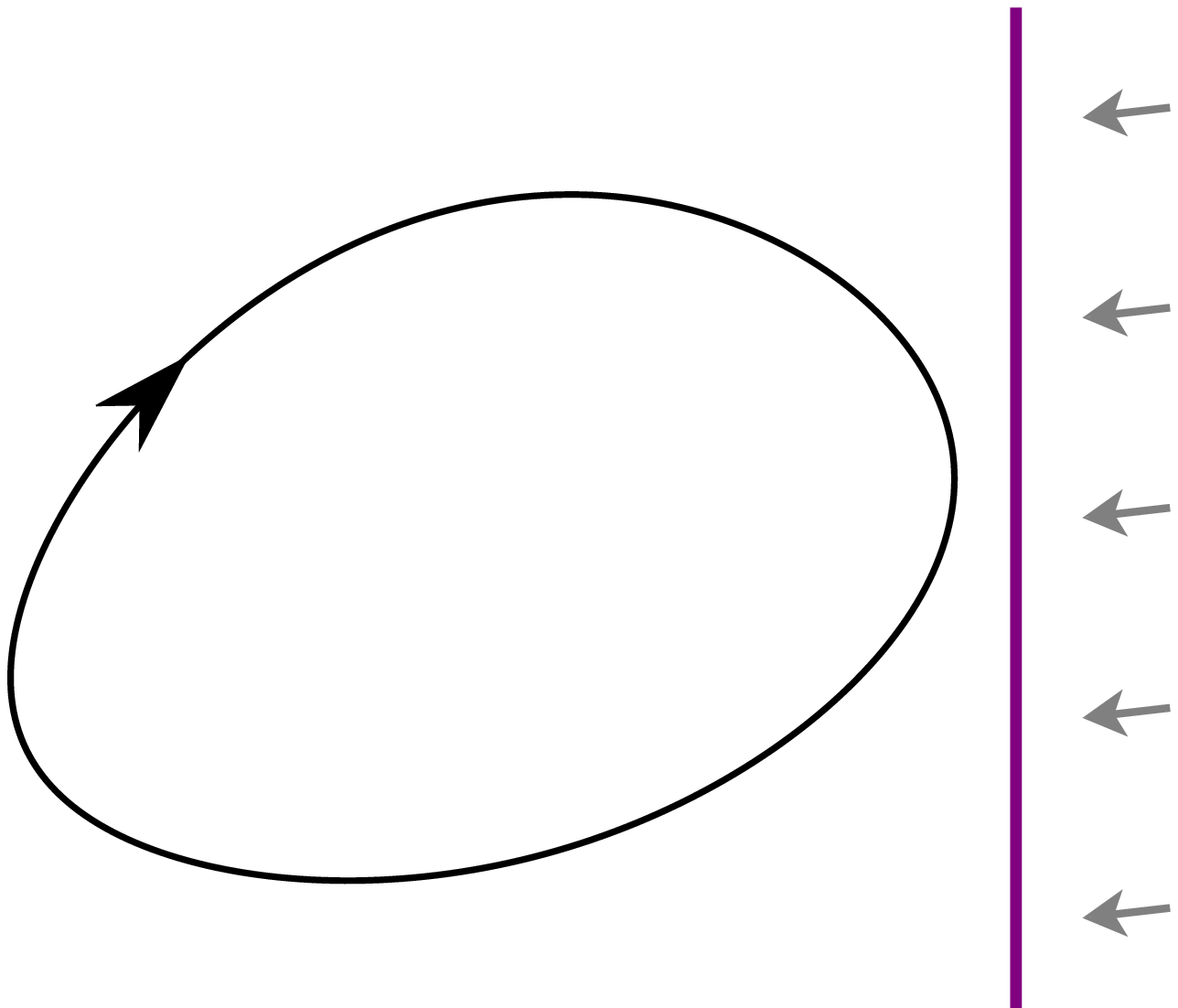}}
\put(5.7,0){\includegraphics[height=4.2cm]{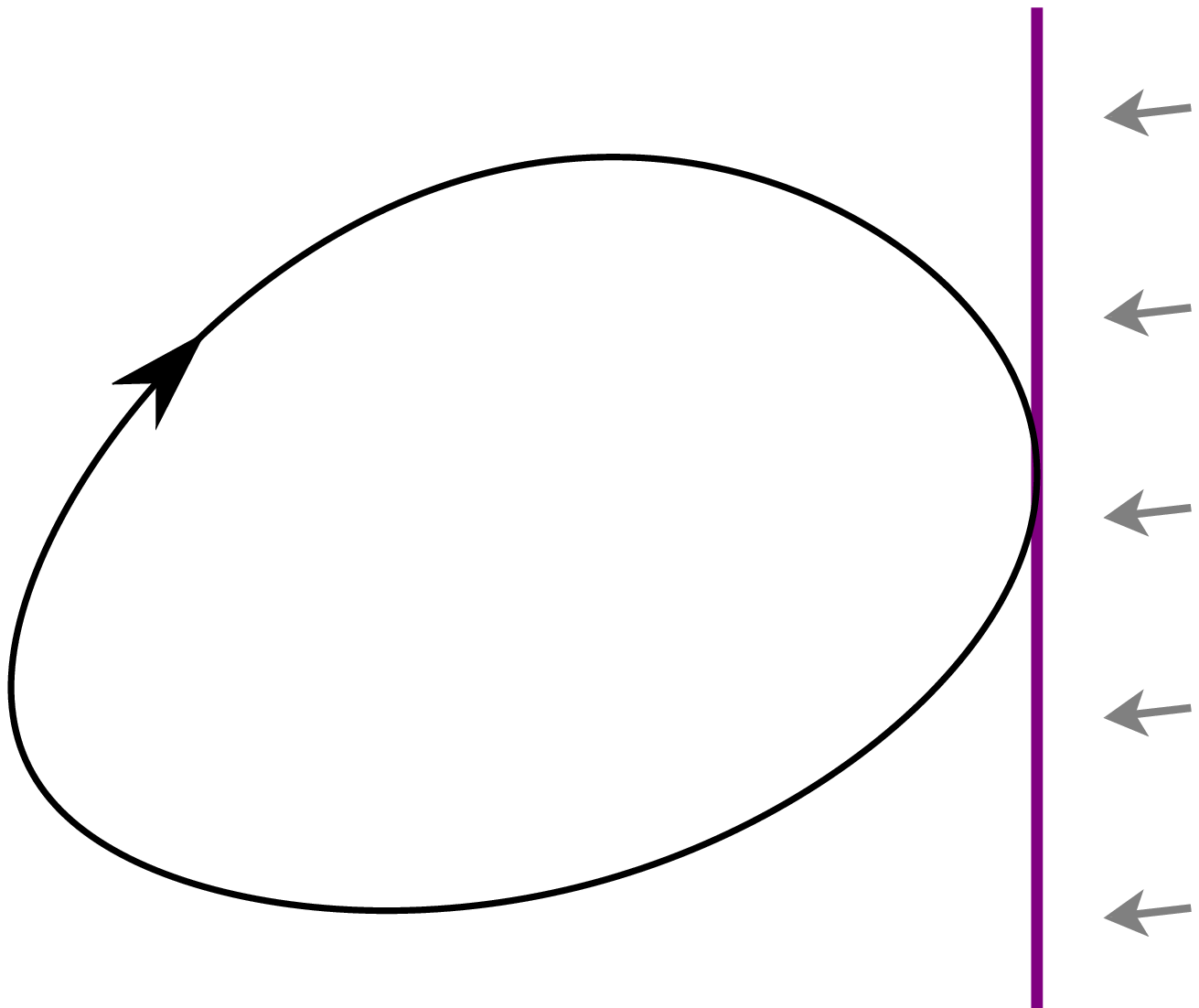}}
\put(11.4,0){\includegraphics[height=4.2cm]{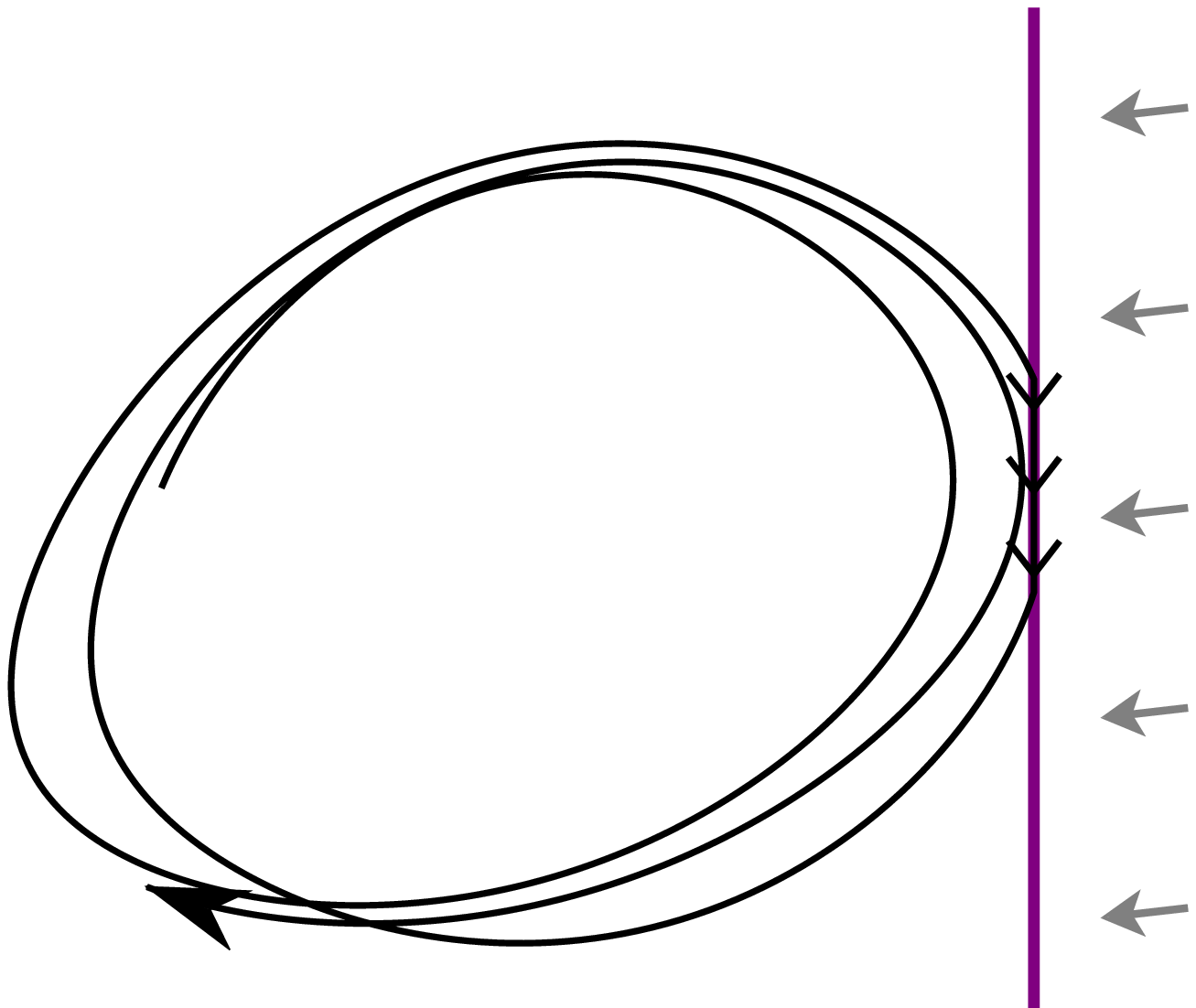}}
\put(1.9,4.1){\small $\gamma < \gamma_{\rm graz}$}
\put(7.6,4.1){\small $\gamma = \gamma_{\rm graz}$}
\put(13.3,4.1){\small $\gamma > \gamma_{\rm graz}$}
\end{picture}
\caption{
Three phase portraits illustrating a grazing-sliding bifurcation
occurring at $\gamma = \gamma_{\rm graz}$, where $\gamma$ is a system parameter.
To the right of a discontinuity surface,
shown with a vertical line, the vector field is directed towards the discontinuity surface, as indicated.
For $\gamma > \gamma_{\rm graz}$ the orbit shown has a segment of sliding motion.
\label{fig:schemGrazSlid}
}
\end{center}
\end{figure}

Grazing-sliding bifurcations occur for piecewise-smooth systems of ODEs
that are discontinuous on manifolds where they are nonsmooth, termed discontinuity surfaces.
At places on discontinuity surfaces where the vector field is directed towards the surface from sides,
orbits evolve on the discontinuity surface --- this is known as sliding motion \cite{DiBu08,Fi88}.
A periodic orbit of a piecewise-smooth system undergoes a bifurcation 
when it grazes a discontinuity surface as parameters are varied.
If, at the point of grazing, the part of the vector field that does not govern the periodic orbit
is directed towards the discontinuity surface,
then this is a grazing-sliding bifurcation, see Fig.~\ref{fig:schemGrazSlid}.
Other bifurcations of this nature are detailed in \cite{CoDi12,JeHo11,KoDi06}.

Grazing-sliding bifurcations arise naturally in mechanical systems with stick-slip friction.
In this context the bifurcation occurs most simply when regular oscillations not involving sticking,
transition to irregular oscillations involving recurring phases of sticking (these correspond to segments of sliding motion),
see for instance \cite{CaGi06,DiKo03,GuHo10,KoPi08,LuGe06} and references within \cite{MaLa12}.
Grazing-sliding bifurcations have been identified in predator-prey models
that are piecewise-smooth due to the assumption that predators
are only harvested when they are in sufficiently high numbers \cite{DeGr03,KuRi03}
and in a two-stage population model \cite{TaLi12}.

The dynamics associated with grazing-sliding bifurcations can be simple or extremely complicated.
An asymptotically stable periodic orbit can simply accumulate a segment of sliding motion.
Alternatively it may bifurcate into an asymptotically stable periodic orbit involving several loops
near the original periodic orbit, some of which involve segments of sliding motion \cite{SzOs09}.
The periodic orbit may bifurcate into a chaotic attractor \cite{Ko05}.
Interestingly, there is no restriction on the dimension of this attractor \cite{Gl15b,GlJe15}\removableFootnote{
Other studies heavily involving sliding bifurcations include \cite{DaNo00,KoDi05,NoKo06}.
}.

In \cite{GlKo12} it was shown that at a grazing-sliding bifurcation an asymptotically stable periodic orbit
can bifurcate into multiple attractors.
More recently in \cite{GlKo16} the same authors introduced an abstract ODE system
for which key calculations could be achieved explicitly and provided 
examples for which an asymptotically stable periodic orbit bifurcates into (i) two asymptotically stable periodic orbits,
and (ii) an asymptotically stable periodic orbit and a chaotic attractor.

The purpose of this paper is to show that infinitely many attractors can be created in grazing-sliding bifurcations.
This is achieved by working with a return map that captures the local dynamics.
The return map is piecewise-smooth because
return trajectories that involve a segment of sliding motion produce a different functional form in the map
than return trajectories that do not.
As was first shown in \cite{DiKo02}, the return map is continuous and piecewise-differentiable.
To leading order, the map can written as
\begin{equation}
f(\bx) = \begin{cases}
A_L \bx + b \mu \;, & e_1^{\sf T} \bx \le 0 \;, \\
A_R \bx + b \mu \;, & e_1^{\sf T} \bx \ge 0 \;,
\end{cases}
\label{eq:f}
\end{equation}
where $\bx \in \mathbb{R}^N$ is the state variable and
$\mu \in \mathbb{R}$ is a parameter. 
The $N \times N$ matrices $A_L$ and $A_R$ differ only in their first columns (by continuity) and $b \in \mathbb{R}^N$.
Here, and throughout the paper, $e_j$ denotes the $j^{\rm th}$ standard basis vector of $\mathbb{R}^N$
for $j = 1,\ldots,N$.
The surface $e_1^{\sf T} \bx = 0$, call it $\Sigma$, is the switching manifold of \eqref{eq:f}.

Let us suppose that the right component of \eqref{eq:f} (the part with $e_1^{\sf T} \bx \ge 0$)
corresponds to return trajectories that involve a segment of sliding motion.
Since sliding motion occurs on a codimension-one surface
(namely the discontinuity surface associated with the grazing-sliding bifurcation),
the range of the right component of \eqref{eq:f} must have dimension less than $N$.
That is, $\det(A_R) = 0$.

The periodic orbit associated with the grazing-sliding bifurcation
corresponds to a fixed point of \eqref{eq:f}.
The grazing-sliding bifurcation occurs for $\mu = 0$ when this fixed point collides with $\Sigma$ at $\bx = \b0$.
In the context of \eqref{eq:f}, this is known as a border-collision bifurcation \cite{Si16}.

A mechanism for the creation of infinitely many attractors in border-collision bifurcations
was introduced for two-dimensional maps in \cite{Si14},
and generalised to maps of any number of dimensions in \cite{SiTu17}.
Here it is shown that this mechanism can occur for grazing-sliding bifurcations.
Although the required codimension is relatively high (the bifurcation is codimension-four instead of codimension-one),
about a point in parameter space at which this phenomenon occurs, for any $K \ge 1$
there exists an open set within which the system has at least $K$ attractors.

The remainder of this paper is organised as follows.
In \S\ref{sec:SiTu17} we introduce symbolic notation to characterise periodic solutions of \eqref{eq:f}.
We then state Theorem \ref{th:SiTu17}, due to \cite{SiTu17}, 
that lists conditions sufficient for \eqref{eq:f} to have infinitely many periodic solutions of a simple type.

It seems that the conditions of Theorem \ref{th:SiTu17} can only be satisfied for \eqref{eq:f} with $\det(A_R) = 0$
if \eqref{eq:f} is at least three-dimensional.
For this reason we focus on \eqref{eq:f} in three dimensions.
In \S\ref{sec:derivingExamples} we describe a practical approach for
determining the parameters of the three-dimensional border-collision normal form
for which the conditions of Theorem \ref{th:SiTu17} may be satisfied.
We then use this approach to construct a two-parameter family of maps satisfying the conditions of Theorem \ref{th:SiTu17}
and numerically obtain two additional examples.

In the spirit of \cite{GlKo16},
we introduce an abstract ODE system that exhibits grazing-sliding bifurcations in \S\ref{sec:odeExample}.
This system is sufficiently simple that
the parameters of the corresponding border-collision normal form
can be written explicitly in terms of the parameters of the ODE system.
Moreover, the system is designed so that
the inverse problem of determining the parameters of the ODE system
that give desired parameters in the normal form can be solved analytically.
This is explained in \S\ref{sec:parameters}
and enables us to generate grazing-sliding bifurcations at which
infinitely many asymptotically stable periodic orbits are generated.
The identification of this phenomenon in mathematical models of real world systems is left for future work.
In \S\ref{sec:bifDiag} we describe the bifurcation diagram for a representative example.
Concluding comments are presented in \S\ref{sec:conc}.

\section{Sufficient conditions for infinitely many attractors}
\label{sec:SiTu17}
\setcounter{equation}{0}

We begin by explaining how periodic solutions of \eqref{eq:f}
can be represented symbolically, as in \cite{Si16}.

Let $\cX \in \{ L,R \}^n$ be a word of length $n$ involving the symbols $L$ and $R$.
We index the elements of such a word from $0$ to $n-1$ and write $\cX = \cX_0 \cdots \cX_{n-1}$.
Given $\cX \in \{ L,R \}^n$ and $\cY \in \{ L,R \}^p$,
the concatenation of $\cX$ and $\cY$ is
\begin{equation}
\cX \cY = \cX_0 \cdots \cX_{n-1} \cY_0 \cdots \cY_{p-1} \;,
\nonumber
\end{equation}
which is a word of length $n+p$.
We write $\cX^k$, where $k \ge 0$ is an integer,
to denote the concatenation of $\cX$ with itself $k$ times. 
For any $i = 0,\ldots,n-1$,
we write $\cX^{\overline{i}}$ to denote the word of length $n$
that equals $\cX$ in all elements except $\cX_i$
(e.g.~if $\cX = RLR$, then $\cX^{\overline{2}} = RLL$).

Let
\begin{align*}
f_L(\bx) &= A_L \bx + b \mu \;, \\
f_R(\bx) &= A_R \bx + b \mu \;,
\end{align*}
denote the two components of $f$, \eqref{eq:f}.
For any $\cX \in \{ L,R \}^n$, let
\begin{equation}
f_\cX = f_{\cX_{n-1}} \circ \cdots \circ f_{\cX_0} \;,
\label{eq:fX}
\end{equation}
denote the composition of $f_L$ and $f_R$ in the order specified by the elements of $\cX$.
The function $f_\cX$ is affine and given by
\begin{equation}
f_\cX(\bx) = M_\cX \bx + P_\cX b \mu \;,
\label{eq:fX2}
\end{equation}
where
\begin{align}
M_\cX &= A_{\cX_{n-1}} \cdots A_{\cX_0} \;, \label{eq:MX} \\
P_\cX &= I + \sum_{i=1}^{n-1} A_{\cX_{n-1}} \cdots A_{\cX_i} \;. \label{eq:PX}
\end{align}
An $n$-tuple $\left( \bx^\cX_0,\ldots,\bx^\cX_{n-1} \right)$ for which
\begin{equation}
f_{\cX_i} \left( \bx^\cX_i \right) = \bx^\cX_{(i+1) \,{\rm mod}\, n} \;, \qquad
{\rm for~all~} i = 0,\ldots,n-1 \;, 
\label{eq:Xcycle}
\end{equation}
is called an $\cX$-cycle.
The $\cX$-cycle is a periodic orbit of $f$ and said to {\em admissible}
if each $\bx^\cX_i$ lies on the ``correct'' side of the switching manifold $\Sigma$, or on $\Sigma$.
To be more precise, for any $\bx \notin \Sigma$ let
\begin{equation}
s(\bx) = \begin{cases}
L \;, & e_1^{\sf T} \bx < 0 \;, \\
R \;, & e_1^{\sf T} \bx > 0 \;.
\end{cases}
\label{eq:s}
\end{equation}
Then the $\cX$-cycle is admissible
if $s \left( \bx^\cX_i \right) = \cX_i$ for all $i$ for which $\bx^\cX_i \notin \Sigma$.
Since $\bx^\cX_0$ is a fixed point of \eqref{eq:fX2},
if no points of an admissible $\cX$-cycle lie on $\Sigma$
then the $\cX$-cycle is asymptotically stable if and only if all eigenvalues of $M_\cX$ have modulus less than $1$.

Given two words $\cX$ and $\cY$, the following result, taken from \cite{SiTu17},
provides sufficient conditions for $f$
to have infinitely many admissible, asymptotically stable $\cX^k \cY$-cycles.

\begin{theorem}
Let $\cX \in \{ L,R \}^n$ and $\cY \in \{ L,R \}^p$
be such that $\cX \cY = \left( \cY \cX \right)^{\overline{0} \overline{\alpha}}$
for some $\alpha \in \{ 1,\ldots,n+p-1 \}$.
\begin{enumerate}
\item
Suppose $M_\cX$ has multiplicity-one eigenvalues $\lambda_1 > 1$ and $\lambda_2 = \frac{1}{\lambda_1}$
and all other eigenvalues of $M_\cX$ have modulus less than $\lambda_2$.
\item
For $j = 1,2$, let $\omega_j^{\sf T}$ and $\zeta_j$ be left and right eigenvectors of $M_\cX$
corresponding to $\lambda_j$ and satisfying $\omega_j^{\sf T} \zeta_j = 1$
(which can always be achieved).
Suppose $e_1^{\sf T} \zeta_1 \ne 0$ and $\lambda_2 < \det(C) < 1$ where
\begin{equation}
C = \begin{bmatrix} \omega_1^{\sf T} \\ \omega_2^{\sf T} \end{bmatrix} M_\cY
\begin{bmatrix} \zeta_1 & \zeta_2 \end{bmatrix},
\label{eq:C}
\end{equation}
is a $2 \times 2$ matrix.
\item
Suppose that the $\cX$-cycle (which must be unique) is an admissible
periodic solution of $f$ with no points on $\Sigma$.
\item
Let $\cS = \cX^\infty \cY \cX^\infty$ be a bi-infinite symbol sequence, with $\cS_0$ corresponding to $\cY_0$.
Suppose there exists an orbit $\{ \by_i \}$ of $f$ that is homoclinic to the $\cX$-cycle and
\begin{enumerate}
\item
$s(\by_i) = \cS_i$ for all $i \in \mathbb{Z}$ for which $\by_i \notin \Sigma$;
\item
$\by_0 = \bx^\cX_0 - \frac{e_1^{\sf T} \bx^\cX_0}{e_1^{\sf T} \zeta_1} \,\zeta_1 \in \Sigma$;
\item
$\by_\alpha \in \Sigma$;
\item
there does not exist $i \ge 0$ for which $\by_i \in \Sigma$ and $\by_{i+n} \in \Sigma$.
\end{enumerate}
\end{enumerate}
Then there exists $k_{\rm min} \ge 0$ such that $f$ has an
admissible, asymptotically stable $\cX^k \cY$-cycle with no points on $\Sigma$
for all $k \ge k_{\rm min}$.
\label{th:SiTu17}
\end{theorem}


\section{The three-dimensional border-collision normal form}
\label{sec:derivingExamples}
\setcounter{equation}{0}



Given a three-dimensional map $f$ of the form \eqref{eq:f}, let
\begin{equation}
\cO_L = \begin{bmatrix}
e_1^{\sf T} A_L^2 \\
e_1^{\sf T} A_L \\
e_1^{\sf T}
\end{bmatrix},
\label{eq:OL}
\end{equation}
and let $\varrho^{\sf T} = e_1^{\sf T} {\rm adj}(I - A_L)$,
where ${\rm adj}(\cdot)$ denotes the {\em adjugate} of a matrix.
If $\det(\cO_L) \ne 0$ (this is the ``observability condition'')
then $f$ can be transformed such that $A_L$, $A_R$, and $b$ have the form
\begin{equation}
\begin{split}
A_L &= \begin{bmatrix} \tau_L & 1 & 0 \\ -\sigma_L & 0 & 1 \\ \delta_L & 0 & 0 \end{bmatrix}, \\
A_R &= \begin{bmatrix} \tau_R & 1 & 0 \\ -\sigma_R & 0 & 1 \\ \delta_R & 0 & 0 \end{bmatrix}, \\
b &= e_1 \;.
\end{split}
\label{eq:bcNormalForm}
\end{equation}
If also $\varrho^{\sf T} b \ne 0$ (a non-degeneracy condition for the vector $b$ in the original map)
then $f$ is conjugate to its transformed version for $\mu \ne 0$ \cite{Si16}.

With \eqref{eq:bcNormalForm} the map \eqref{eq:f}
is known the three-dimensional border-collision normal form.
The parameters $\tau_{L,R}$, $\sigma_{L,R}$, and $\delta_{L,R}$ are conveniently the
trace, second trace, and determinant of $A_{L,R}$, see Appendix \ref{app:secondTrace}.

In this section we work with the three-dimensional border-collision normal form.
We fix $\delta_R = 0$, so that $\det(A_R) = 0$, 
and search for values of $\tau_L, \tau_R, \sigma_L, \sigma_R, \delta_L, \in \mathbb{R}$
that satisfy the conditions of Theorem \ref{th:SiTu17} for $\mu = 1$ and some words $\cX$ and $\cY$.


\subsection{Determining parameter values that give infinitely many attractors}
\label{sub:generalProcedure}

The phenomenon described by Theorem \ref{th:SiTu17} is codimension-three
because $\lambda_2 = \frac{1}{\lambda_1}$,
$\by_\alpha \in \Sigma$,
and the requirement that $\by_0$ belongs to the stable manifold of the $\cX$-cycle,
are independent codimension-one conditions.
It is not particularly helpful to {\em directly} use these conditions
to generate restrictions on the values of $\tau_L$, $\tau_R$, $\sigma_L$, $\sigma_R$, and $\delta_L$,
because, for instance, the eigenvalues $\lambda_1$ and $\lambda_2$
are given by the roots of a quadratic equation (assuming $\delta_R = 0$)
and the resulting square-roots create expressions that seem to be too complicated to deal with.

Instead we derive three alternate conditions
that lead to polynomial restrictions on the parameter values.
This was done for the two-dimensional border-collision normal form in \cite{Si14}.
Here we state merely state these conditions; their derivation in given in Appendix \ref{app:construction}.
They are not intended to provide additional insight into the phenomenon described by Theorem \ref{th:SiTu17},
only to be used as a tool for finding suitable parameter values.
Indeed their solutions may not satisfy all conditions of Theorem \ref{th:SiTu17}.

Let the words $\cX$ and $\cY$ be given, where $\cX \cY = \left( \cY \cX \right)^{\overline{0} \overline{\alpha}}$
for some $\alpha \in \{ 1,\ldots,n+p-1 \}$.
Here we assume $\cX$ and $\cY$ both end with the symbol $R$,
as this provides useful simplification but is not too restrictive.
We can first use the conditions of Theorem \ref{th:SiTu17} to calculate the point $\by_0$.
We have $\by_0 \in \Sigma$ (by definition), thus the first component of $\by_0$ is zero.
We have $\by_0 = f_R(\by_{-1})$ (because $\cX$ ends in $R$),
thus the third component of $\by_0$ is zero (because $A_R$ is given by \eqref{eq:bcNormalForm} with $\delta_R = 0$).
Finally, the second component of $\by_0$ can be determined by the condition $\by_\alpha = f^\alpha(\by_0) \in \Sigma$.
Specifically, from \eqref{eq:fX2} we obtain
\begin{equation}
e_2^{\sf T} \by_0 = \frac{-e_1^{\sf T} P_{\tilde{\cX}} b \mu}{e_1^{\sf T} M_{\tilde{\cX}} e_2} \;,
\label{eq:y02}
\end{equation}
where $\tilde{\cX}$ denotes the first $\alpha$ elements of $\cX \cY$.

Once $\by_0$ is calculated, let
\begin{align}
\psi_1 &= P_\cX e_1 \mu - (I - M_\cX) \by_0 \;, \label{eq:psi1} \\
\psi_2 &= M_\cY \psi_1 \;, \label{eq:psi2}
\end{align}
and
\begin{align}
\xi_1 &= M_\cX \psi_1 e_1^{\sf T} \psi_1 - \psi_1 e_1^{\sf T} M_\cX \psi_1 \;, \label{eq:xi1} \\
\xi_2 &= M_\cX \psi_2 e_1^{\sf T} M_\cX \psi_1 - \psi_2 e_1^{\sf T} \psi_1 \;. \label{eq:xi2}
\end{align}
Then our three alternate conditions are
\begin{align}
\sigma_\cX &= 1 \;, \label{eq:construct1} \\
e_2^{\sf T} \xi_1 &= 0 \;, \label{eq:construct2} \\
e_1^{\sf T} \xi_2 &= 0 \;, \label{eq:construct3}
\end{align}
where $\sigma_\cX$ denotes the second trace of $M_\cX$.
Instances of the denominator of \eqref{eq:y02} that arise in
\eqref{eq:construct1}--\eqref{eq:construct3} can be factored out
leaving equations that are polynomial in $\tau_L$, $\tau_R$, $\sigma_L$, $\sigma_R$, and $\delta_L$.

In summary, in order to find values of the parameters in \eqref{eq:bcNormalForm}
for which $f$ has infinitely many admissible, asymptotically stable $\cX^k \cY$-cycles,
we solve \eqref{eq:construct1}-\eqref{eq:construct3} (derived in Appendix \ref{app:construction})
and check that all conditions of Theorem \ref{th:SiTu17} are satisfied.

\subsection{Calculations with $\cX = RLR$ and $\cY = LR$}
\label{sub:derivationMainExample}

Here we consider
\begin{equation}
\cX = RLR \;, \qquad \cY = LR \;,
\label{eq:XYF}
\end{equation}
for which $\cX \cY = \left( \cY \cX \right)^{\overline{0} \overline{\alpha}}$ for $\alpha = 1$.
With $\mu = 1$, for this value of $\alpha$ we have $\by_0 = [0,-1,0]^{\sf T}$.
With $\delta_R = 0$, the second trace of $M_\cX = A_R A_L A_R$
is $\sigma_{\cX} = (\sigma_L \sigma_R - \delta_L \tau_R) \sigma_R$.
Thus \eqref{eq:construct1} gives
\begin{equation}
\delta_L = \frac{\sigma_L \sigma_R^2 - 1}{\tau_R \sigma_R} \;.
\label{eq:dLex1}
\end{equation}
Using a symbolic toolbox, numerically we found that
$e_2^{\sf T} \xi_1$ is an affine function of $\tau_L$ (for this example)
and \eqref{eq:construct2} can be rearranged to produce
\begin{align}
\tau_L &=
\frac{1}{\sigma_R^2 - \sigma_R \tau_R^2 - \sigma_R - \tau_R^3 - \tau_R^2 - \tau_R}
\Big( \delta_L \sigma_R - \sigma_L - \delta_L - 2 \delta_L \tau_R
+ 2 \sigma_L \sigma_R - \sigma_L \tau_R - \sigma_R \tau_R \nonumber \\
&\quad- 2 \delta_L \tau_R^2 - \delta_L \tau_R^3 - \sigma_L \sigma_R^2 - \sigma_L \tau_R^2 - \sigma_R \tau_R^2 - \sigma_R^2 \tau_R - \sigma_R^2 + \sigma_L \sigma_R \tau_R^2 + \delta_L \sigma_R \tau_R + \sigma_L \sigma_R \tau_R \Big) \;.
\label{eq:tLex1}
\end{align}
The quantity $e_1^{\sf T} \xi_2$ contains too many terms to be given here,
but upon substituting \eqref{eq:dLex1} and \eqref{eq:tLex1} simplifies to 
\begin{equation}
e_1^{\sf T} \xi_2 = 
\frac{\tau_R (\tau_R + 1)^2 (1 - \sigma_R)
(\sigma_R^2 + \sigma_R \tau_R - \sigma_R + \tau_R^2 + \tau_R + 1)
(\sigma_R - \tau_R^2 - \tau_R - 1)^4 (\tau_R + \sigma_R + 1)}
{\sigma_R^2 (\sigma_R^2 - \sigma_R \tau_R^2 - \sigma_R - \tau_R^3 - \tau_R^2 - \tau_R)^3} \;.
\label{eq:xi2ex1}
\end{equation}
In view of \eqref{eq:construct3}, we require one factor in the numerator of \eqref{eq:xi2ex1} to be zero.
By considering each factor in turn we find that
all conditions of Theorem \ref{th:SiTu17} can only be satisfied if the last factor is zero, that is
\begin{equation}
\tau_R = -(\sigma_R + 1) \;.
\label{eq:tRex1}
\end{equation}
Then by substituting \eqref{eq:tRex1} into \eqref{eq:dLex1} we obtain
\begin{equation}
\delta_L = \frac{1 - \sigma_L \sigma_R^2}{\sigma_R (\sigma_R + 1)} \;.
\label{eq:dLex1b}
\end{equation}
Lastly by substituting \eqref{eq:tRex1} and \eqref{eq:dLex1b} into \eqref{eq:tLex1} we obtain
\begin{equation}
\tau_L = \frac{1}{\sigma_R^2+1} - \frac{\sigma_L+\sigma_R}{\sigma_R+1} \;.
\label{eq:tLex1b}
\end{equation}

\subsection{A two-parameter family}
\label{sub:verificationMainExample}

\begin{figure}[b!]
\begin{center}
\setlength{\unitlength}{1cm}
\begin{picture}(8,6)
\put(0,0){\includegraphics[height=6cm]{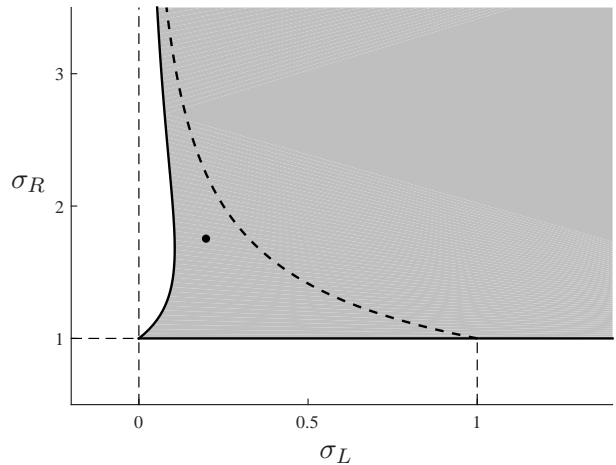}}
\put(4.1,0){\small $\sigma_L$}
\put(0,3.6){\small $\sigma_R$}
\end{picture}
\caption{
The domain $\cD$ \eqref{eq:D}.
The dashed curve, $\sigma_L = \frac{1}{\sigma_R^2}$, is where $\delta_L = 0$.
The black dot is the point $(\sigma_L,\sigma_R) = (0.2,1.75)$
used in Figs.~\ref{fig:qqExdRzero6} and \ref{fig:bifDiagGrazSliding}.
\label{fig:domainExF}
}
\end{center}
\end{figure}

To satisfy the conditions of Theorem \ref{th:SiTu17} with $\cX = RLR$ and $\cY = LR$,
equations \eqref{eq:tRex1}--\eqref{eq:tLex1b} must hold.
Here we show that the conditions of Theorem \ref{th:SiTu17} are indeed satisfied
if the values of the two undetermined parameters, $\sigma_L$ and $\sigma_R$,
belong to the domain
\begin{equation}
\cD = \left\{ (\sigma_L,\sigma_R) ~\middle|~
\sigma_L > \frac{\sigma_R - 1}{\sigma_R \left( \sigma_R^2 + 1 \right)} ,\,
\sigma_R > 1 \right\},
\label{eq:D}
\end{equation}
shown in Fig.~\ref{fig:domainExF}.

\begin{theorem}
Choose any $(\sigma_L,\sigma_R) \in \cD$,
let $\tau_R$, $\delta_L$, and $\tau_L$ be given by \eqref{eq:tRex1}--\eqref{eq:tLex1b},
and let $\delta_R = 0$ and $\mu = 1$.
Then there exists $k_{\rm min} \in \mathbb{Z}$ such that
for $\cX = RLR$ and $\cY = LR$ the map \eqref{eq:f} with \eqref{eq:bcNormalForm}
has an admissible, asymptotically stable $\cX^k \cY$-cycle
with no points on $\Sigma$ for all $k \ge k_{\rm min}$.
\label{th:twoParamEx}
\end{theorem}

Theorem \ref{th:twoParamEx} is proved in Appendix \ref{app:proof}
by simply showing that all conditions of Theorem \ref{th:SiTu17} are satisfied.
Theorem \ref{th:twoParamEx} can also be proved by calculating the $\cX^k \cY$-cycles directly, as in \cite{Si14}.
The latter approach requires lengthy calculations, and so is not included here,
but reveals that we can take $k_{\rm min} = 1$ for all $(\sigma_L,\sigma_R) \in \cD$.

\begin{figure}[b!]
\begin{center}
\setlength{\unitlength}{1cm}
\begin{picture}(16,8)
\put(0,0){\includegraphics[height=8cm]{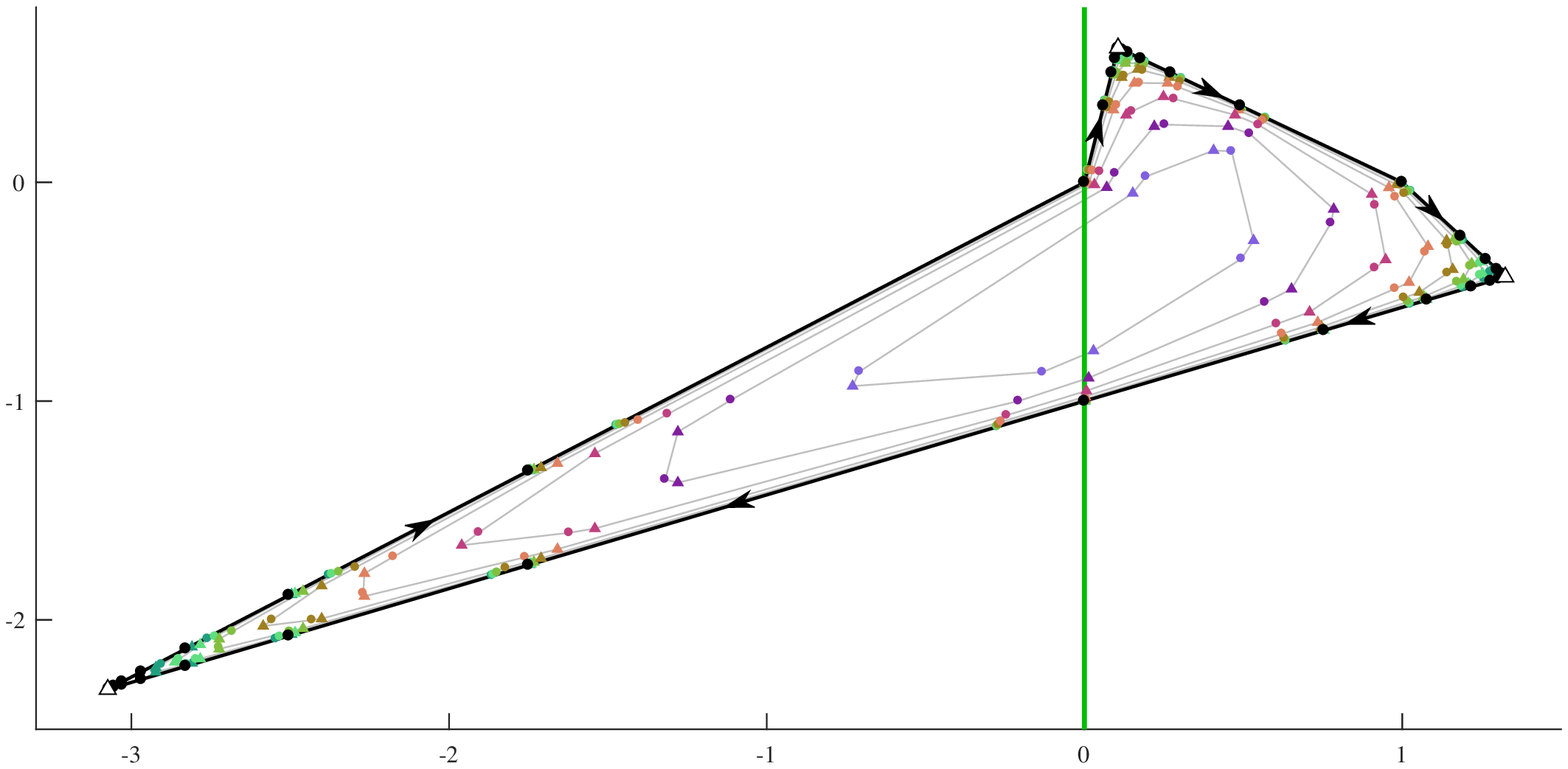}}
\put(9.45,0){\small $e_1^{\sf T} \bx$}
\put(0,5.1){\small $e_2^{\sf T} \bx$}
\put(15.48,5.1){\footnotesize $\bx^{\cX}_0$}
\put(1.06,1.3){\footnotesize $\bx^{\cX}_1$}
\put(11.66,7.6){\footnotesize $\bx^{\cX}_2$}
\put(13.46,4.53){\scriptsize $\by_{\hspace{.8mm}3}$}
\put(13.64,4.508){\tiny -}
\put(5.4,3.53){\scriptsize $\by_{\hspace{.8mm}2}$}
\put(5.58,3.508){\tiny -}
\put(12.62,7.16){\scriptsize $\by_{\hspace{.8mm}1}$}
\put(12.80,7.138){\tiny -}
\put(11.29,3.88){\scriptsize $\by_0$}
\put(10.8,6.3){\scriptsize $\by_1$}
\put(14.42,6.3){\scriptsize $\by_2$}
\put(5.68,2.22){\scriptsize $\by_3$}
\put(11.29,1.7){\small $\Sigma$}
\end{picture}
\caption{
A phase portrait of \eqref{eq:f} with \eqref{eq:bcNormalForm}, \eqref{eq:exF}, and $\mu = 1$.
Here $\cX = RLR$ and $\cY = LR$.
Asymptotically stable $\cX^k \cY$-cycles [saddle-type $\cX^k \cY^{\overline{0}}$-cycles]
are shown with coloured circles [triangles] for $k = 1,\ldots,8$.
The saddle-type $\cX$-cycle is shown with unshaded triangles.
\label{fig:qqExdRzero6}
}
\end{center}
\end{figure}

As a specific example, consider the values $(\sigma_L,\sigma_R) = (0.2,1.75)$.
From \eqref{eq:tRex1}--\eqref{eq:tLex1b}, altogether we have
\begin{equation}
\begin{aligned}
\tau_L &= -\frac{331}{715} \;, &
\tau_R &= -\frac{11}{4} \;, \\
\sigma_L &= \frac{1}{5} \;, &
\sigma_R &= \frac{7}{4} \;, \\
\delta_L &= \frac{31}{385} \;, &
\delta_R &= 0 \;.
\end{aligned}
\label{eq:exF}
\end{equation}
Fig.~\ref{fig:qqExdRzero6} shows a phase portrait using these values.
This figure shows the $\cX^k \cY$-cycles for $k = 1,\ldots,8$ (with circles).
For these values of $k$, saddle-type $\cX^k \cY^{\overline{0}}$-cycles also exist (shown with triangles).
It seems typical for the stable manifolds of these saddle solutions
to form the boundaries of the basins of attraction of the $\cX^k \cY$-cycles, see \cite{Si14}.
To show the $\cX^k \cY$ and $\cX^k \cY^{\overline{0}}$-cycles clearly,
in Fig.~\ref{fig:qqExdRzero6} for each $k$ the points of these periodic solutions are connected by line segments.

The $\cX$-cycle (with points $\bx^\cX_0$, $\bx^\cX_1$, and $\bx^\cX_2$)
has a one-dimensional unstable manifold and a two-dimensional stable manifold.
As shown in \cite{SiTu17}, the branch of the unstable manifold of the $\cX$-cycle that contains the homoclinic orbit $\{ \by_i \}$
is a subset of the stable manifold of the $\cX$-cycle.
This branch is indicated with solid black lines in Fig.~\ref{fig:qqExdRzero6}.
There also exists an asymptotically stable $\cX^{\overline{2}}=RLL$-cycle,
but this is not visible in Fig.~\ref{fig:qqExdRzero6} as it lies outside the region of phase space shown.

\subsection{Additional examples}
\label{sub:additionalExamples}

\begin{figure}[b!]
\begin{center}
\setlength{\unitlength}{1cm}
\begin{picture}(16,8)
\put(0,0){\includegraphics[height=8cm]{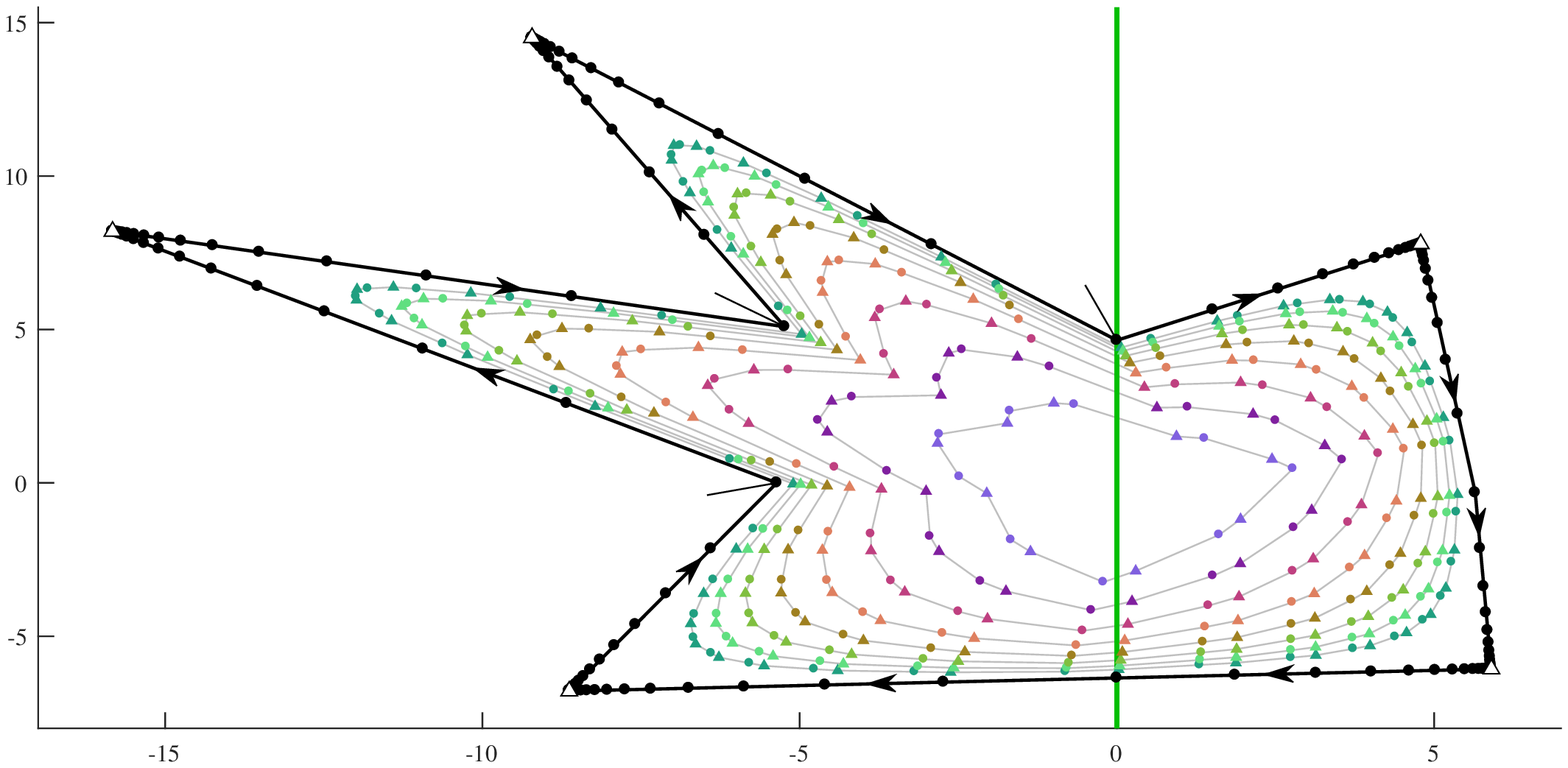}}
\put(9.45,0){\small $e_1^{\sf T} \bx$}
\put(0,5.1){\small $e_2^{\sf T} \bx$}
\put(15.36,1.26){\footnotesize $\bx^{\cX}_0$}
\put(5.53,1.03){\footnotesize $\bx^{\cX}_1$}
\put(1.33,5.95){\footnotesize $\bx^{\cX}_2$}
\put(5.2,7.65){\footnotesize $\bx^{\cX}_3$}
\put(14.64,5.62){\footnotesize $\bx^{\cX}_4$}
\put(11.59,1.09){\scriptsize $\by_0$}
\put(7.1,3.06){\scriptsize $\by_1$}
\put(7.16,5.14){\scriptsize $\by_2$}
\put(10.99,5.31){\scriptsize $\by_3$}
\put(15.15,3.11){\scriptsize $\by_4$}
\put(11.59,7.5){\small $\Sigma$}
\end{picture}
\caption{
A phase portrait of \eqref{eq:f} with \eqref{eq:bcNormalForm}, \eqref{eq:ex20}, and $\mu = 1$
using the same conventions as Fig.~\ref{fig:qqExdRzero6}.
The $\cX^k \cY$-cycles and $\cX^k \cY^{\overline{0}}$-cycles
are shown for $k = 1,\ldots,8$, where $\cX$ and $\cY$ are given by \eqref{eq:XYex20}.
\label{fig:qqExdRzero20}
}
\end{center}
\end{figure}

Here we provide two numerical examples
using combinations of words $\cX$ and $\cY$ that have not been treated
in previous studies of this phenomenon.

With
\begin{equation}
\cX = RLLLR \;, \qquad
\cY = LLLR \;,
\label{eq:XYex20}
\end{equation}
for which $\alpha = 3$, we fixed the values of $\sigma_L$ and $\sigma_R$
and solved \eqref{eq:construct1}--\eqref{eq:construct3} numerically to obtain
\begin{equation}
\begin{aligned}
\tau_L &= 1.1634777991 \;, &
\tau_R &= -0.6037872000 \;, \\
\sigma_L &= 0.95 \;, &
\sigma_R &= 1.15 \;, \\
\delta_L &= 0.0608806824 \;, &
\delta_R &= 0 \;,
\end{aligned}
\label{eq:ex20}
\end{equation}
accurate to ten decimal places.
Fig.~\ref{fig:qqExdRzero20} shows a phase portrait using these values.
Here admissible, asymptotically stable $\cX^k \cY$-cycles exist for at least $k = 1,\ldots,8$.
We expect that with the exact solution to \eqref{eq:construct1}--\eqref{eq:construct3},
all conditions of Theorem \ref{th:SiTu17} are satisfied
and thus infinitely many admissible, asymptotically stable $\cX^k \cY$-cycles exist.
The values of $\sigma_L$ and $\sigma_R$ in \eqref{eq:ex20} were obtained via numerical exploration.

\begin{figure}[t!]
\begin{center}
\setlength{\unitlength}{1cm}
\begin{picture}(16,8)
\put(0,0){\includegraphics[height=8cm]{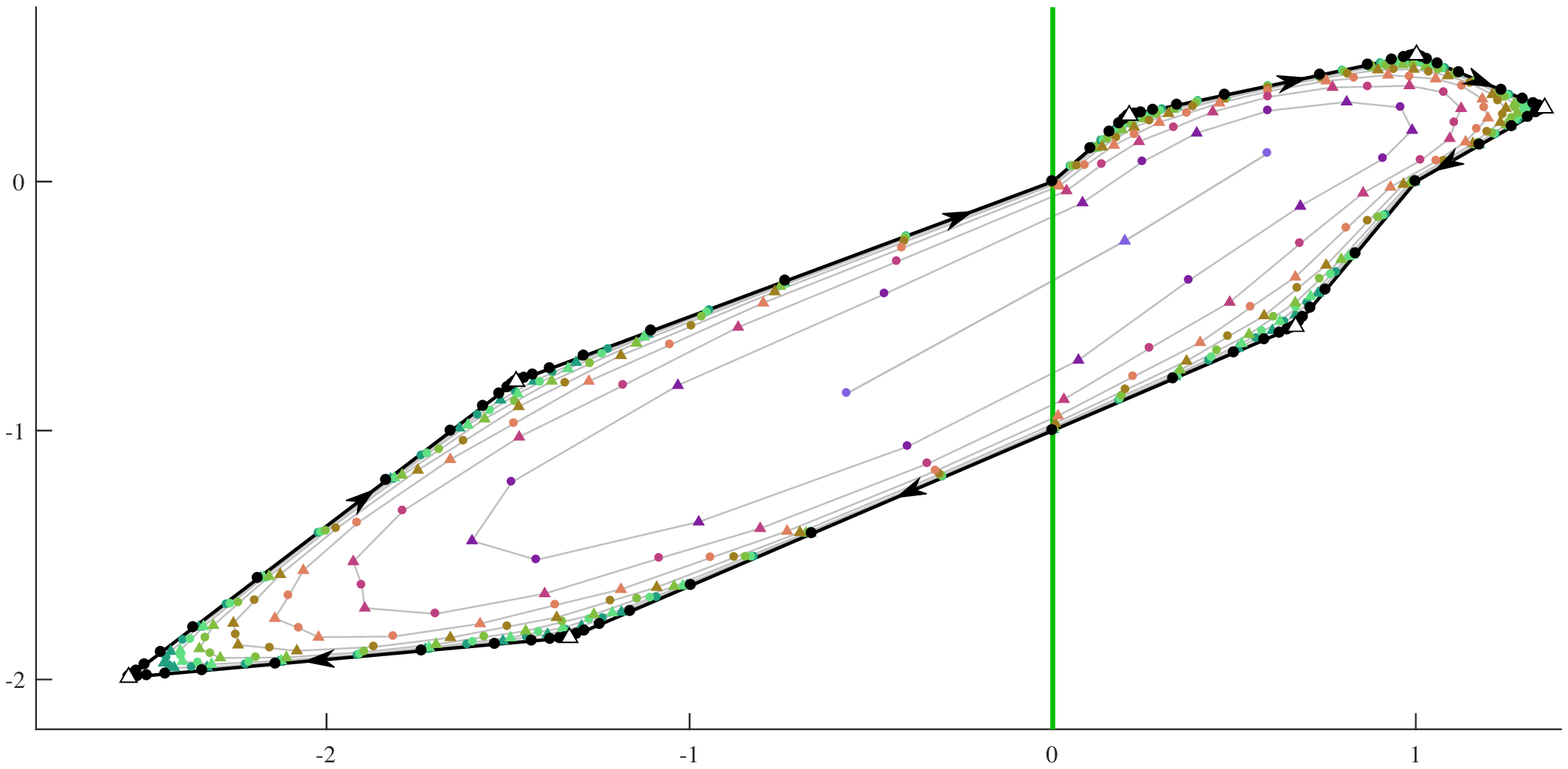}}
\put(9.45,0){\small $e_1^{\sf T} \bx$}
\put(0,5.1){\small $e_2^{\sf T} \bx$}
\put(13.42,4.67){\footnotesize $\bx^{\cX}_0$}
\put(5.04,4.32){\footnotesize $\bx^{\cX}_1$}
\put(15.54,6.55){\footnotesize $\bx^{\cX}_2$}
\put(1.24,1.35){\footnotesize $\bx^{\cX}_3$}
\put(14.35,7.65){\footnotesize $\bx^{\cX}_4$}
\put(6,1.4){\footnotesize $\bx^{\cX}_5$}
\put(11.44,7.11){\footnotesize $\bx^{\cX}_6$}
\put(10.95,3.58){\scriptsize $\by_0$}
\put(10.51,6.32){\scriptsize $\by_1$}
\put(10.95,1.7){\small $\Sigma$}
\end{picture}
\caption{
A phase portrait of \eqref{eq:f} with \eqref{eq:bcNormalForm}, \eqref{eq:ex22}, and $\mu = 1$
using the same conventions as Fig.~\ref{fig:qqExdRzero6}.
The $\cX^k \cY$-cycles and $\cX^k \cY^{\overline{0}}$-cycles
are shown for $k = 0,\ldots,7$, where $\cX$ and $\cY$ are given by \eqref{eq:XYex22}.
\label{fig:qqExdRzero22}
}
\end{center}
\end{figure}

With
\begin{equation}
\cX = RLRLRLR \;, \qquad
\cY = LR \;,
\label{eq:XYex22}
\end{equation}
for which $\alpha = 1$,
we solved \eqref{eq:construct1}--\eqref{eq:construct3} numerically to obtain
\begin{equation}
\begin{aligned}
\tau_L &= -0.7831707737 \;, &
\tau_R &= -2.8347004550 \;, \\
\sigma_L &= 0.2 \;, &
\sigma_R &= 1.2 \;, \\
\delta_L &= 0.2473051527 \;, &
\delta_R &= 0 \;,
\end{aligned}
\label{eq:ex22}
\end{equation}
accurate to ten decimal places.
Admissible, asymptotically stable $\cX^k \cY$-cycles exist for at least $k = 0,\ldots,7$,
as shown in Fig.~\ref{fig:qqExdRzero22}.
Note that with $k = 0$, we have $\cX^k \cY^{\overline{0}} = RR$.
The $RR$-cycle consists only of the fixed point of $f_R$.

\section{An abstract ODE system}
\label{sec:odeExample}
\setcounter{equation}{0}

Here we study the three-dimensional non-autonomous system
\begin{equation}
\begin{bmatrix} \dot{X} \\ \dot{Y} \\ \dot{Z} \end{bmatrix} =
\begin{cases}
\begin{bmatrix} Y \\ Z \\ -\alpha_1 (X+1) - \alpha_2 Y - \alpha_3 Z + \gamma \cos(t) \end{bmatrix}, & X < 0 \;, \\
\begin{bmatrix} -1 \\ \beta_1 \\ \beta_2 \end{bmatrix}, & X > 0 \;,
\end{cases}
\label{eq:odeEx}
\end{equation}
where $\alpha_1, \alpha_2, \alpha_3, \beta_1, \beta_2, \gamma \in \mathbb{R}$ are constants.
The system \eqref{eq:odeEx} is piecewise-smooth with the discontinuity surface $X = 0$
and we write $\bX = (X,Y,Z)$.
We treat $\gamma$ as the primary bifurcation parameter.
This parameter can be thought of as a forcing amplitude;
indeed \eqref{eq:odeEx} is motivated by a harmonically forced linear oscillator.
We have included five additional parameters
so that we can fit these to five given non-zero eigenvalues of $A_L$ and $A_R$, as achieved in \S\ref{sec:parameters}.
While $X < 0$, the explicit solution to \eqref{eq:odeEx} is available.
This facilitates accurate numerical simulations, presented in \S\ref{sec:bifDiag}.
In this section we use the explicit solution to identify a grazing-sliding bifurcation
and calculate the return map to leading order.

With $\gamma = 0$, the point $\bX = (-1,0,0)$ is an equilibrium of \eqref{eq:odeEx}.
Assuming $(\alpha_1 - \alpha_3)^2 + (\alpha_2 - 1)^2 \ne 0$,
for sufficiently small $\gamma > 0$ the system \eqref{eq:odeEx} has an oscillatory solution
in the left half-space ($X < 0$) centred at $\bX = (-1,0,0)$.
This solution is given by
\begin{equation}
\bX_p(t) = \frac{\gamma}{(\alpha_1 - \alpha_3)^2 + (\alpha_2 - 1)^2}
\left( \begin{bmatrix} \alpha_1 - \alpha_3 \\ \alpha_2 - 1 \\ -(\alpha_1 - \alpha_3) \end{bmatrix} \cos(t) +
\begin{bmatrix} \alpha_2 - 1 \\ -(\alpha_1 - \alpha_3) \\ -(\alpha_2 - 1) \end{bmatrix} \sin(t) \right) - e_1 \;,
\label{eq:Xp}
\end{equation}
and grazes $X = 0$ at
\begin{equation}
\gamma_{\rm graz} = \sqrt{(\alpha_1 - \alpha_3)^2 + (\alpha_2 - 1)^2} \;.
\label{eq:gammagraz}
\end{equation}

In order to employ standard techniques regarding grazing events of piecewise-smooth systems,
we reinterpret \eqref{eq:odeEx} as a four-dimensional autonomous system by treating $t$ as a variable
(i.e.~with $\dot{t} = 1$).
We also take $t$ modulo $2 \pi$,
so that in the cylindrical phase space $\mathbb{R}^3 \times \mathbb{S}$
the oscillatory solution $\bX_p(t)$ is a periodic orbit.
Grazing occurs at the point
\begin{align}
\bX_{\rm graz} &= (0,0,-1), \label{eq:Xgraz} \\
t_{\rm graz} &= \tan^{-1} \left( \frac{\alpha_2 - 1}{\alpha_1 - \alpha_3} \right), \label{eq:tgraz}
\end{align}
where $t_{\rm graz} \in (0,\pi)$ if $\alpha_2 - 1 > 0$ and
$t_{\rm graz} \in (\pi,2 \pi)$ if $\alpha_2 - 1 < 0$.
Note that $\gamma = \gamma_{\rm graz}$ is a grazing-sliding bifurcation
because at the point of grazing the right half-system is directed towards the discontinuity surface
(specifically $\dot{X} = -1$).

For $Y > 0$, orbits slide on the discontinuity surface $X = 0$
because both components of \eqref{eq:odeEx} are directed towards $X = 0$.
As detailed in \cite{DiBu08,Fi88},
this sliding motion is governed by the convex combination of the components of \eqref{eq:odeEx}
that is tangent to $X = 0$:
\begin{equation}
\begin{bmatrix} \dot{Y} \\ \dot{Z} \end{bmatrix} =
\frac{1}{Y + 1}
\begin{bmatrix}
\beta_1 Y + Z \\
-\alpha_1 + (\beta_2 - \alpha_2) Y - \alpha_3 Z + \gamma \cos(t)
\end{bmatrix},
\label{eq:slidingODE}
\end{equation}
and $\dot{t} = 1$.
For $Y < 0$, orbits cross $X = 0$ and enter the left half-space.

\begin{figure}[t!]
\begin{center}
\setlength{\unitlength}{1cm}
\begin{picture}(8,6)
\put(0,0){\includegraphics[height=6cm]{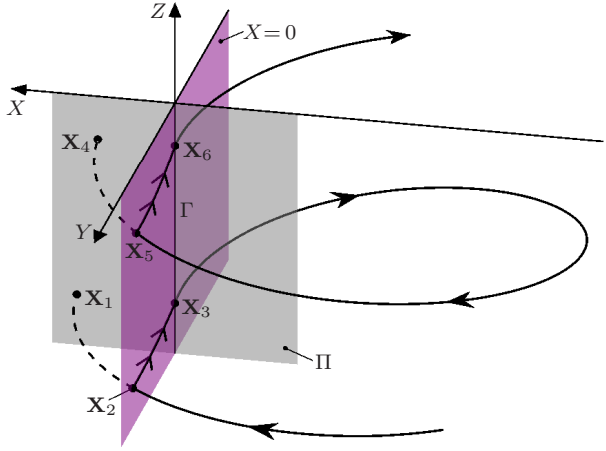}}
\put(0,4.46){\scriptsize $X$}
\put(.92,2.81){\scriptsize $Y$}
\put(1.91,5.75){\scriptsize $Z$}
\put(2.31,3.1){\scriptsize $\Gamma$}
\put(4.11,1.07){\scriptsize $\Pi$}
\put(3.15,5.47){\scriptsize $X \hspace{-.8mm}=\hspace{-.3mm} 0$}
\put(1.04,1.93){\scriptsize $\bX_1$}
\put(1.07,.5){\scriptsize $\bX_2$}
\put(2.33,1.81){\scriptsize $\bX_3$}
\put(.73,4){\scriptsize $\bX_4$}
\put(1.58,2.55){\scriptsize $\bX_5$}
\put(2.35,3.9){\scriptsize $\bX_6$}
\end{picture}
\caption{
Part of a typical orbit of \eqref{eq:odeEx} near the grazing-sliding bifurcation $\gamma = \gamma_{\rm graz}$.
Take care to note the reverse orientation of the $X$-axis.
\label{fig:schemGrazSlidReturnMap}
}
\end{center}
\end{figure}

Let $\Pi$ denote the Poincar\'{e} section $Y = 0$.
Orbits cease sliding and enter $X < 0$ at the intersection of $\Pi$ with $X = 0$, call it $\Gamma$.
Here we use $\Pi$ to build a return map valid near the grazing-sliding bifurcation.
To do this we use the standard approach of combining a global map with a discontinuity map, see \cite{DiBu08,DiKo02,GlKo12}.

Fig.~\ref{fig:schemGrazSlidReturnMap} shows part of a typical orbit near the grazing-sliding bifurcation.
The orbit intersects $X = 0$ at $\bX_2$,
then slides along to $\bX_3 \in \Gamma$.
It then sojourns away from $X = 0$ (following close to the path of $\bX_p(t)$ at $\gamma = \gamma_{\rm graz}$),
intersects $X = 0$ at $\bX_5$, then slides along to $\bX_6 \in \Gamma$.

As shown in Fig.~\ref{fig:schemGrazSlidReturnMap},
we extend the orbit beyond $\bX_2$ and $\bX_5$ to
the virtual points $\bX_1$ and $\bX_4$ where the orbit would intersect $\Pi$
if it were governed by the left half-system of \eqref{eq:odeEx} in $X > 0$.
The {\em global map}, $\cP_g : \Pi \to \Pi$,
is defined as the next intersection of the orbit with $\Pi$ obtained by just using the left half-system.
That is, $\cP_g(\bX_3) = \bX_4$.
The {\em discontinuity map} $\cP_d : \Pi \to \Pi$
is defined as the necessary correction to generate the true point of intersection with $\Pi$.
That is, $\cP_d(\bX_4) = \bX_6$.
For points $\bX \in \Pi$ with $X < 0$, we take $\cP_d$ to be the identity map.

The composition $\cP_d \circ \cP_g$ provides the true return map on $\Pi$.
This form is convenient because $\cP_d$ is a local map and can be computed via asymptotic expansions,
while $\cP_g$ involves transversal intersections with a single Poincar\'{e} section and only one functional form of \eqref{eq:odeEx}.
Below we work with the alternate return map
\begin{equation}
\cP = \cP_g \circ \cP_d \;,
\label{eq:P}
\end{equation}
as this ordering allows for a simpler description of the switching manifold.
The map $\cP$ captures the dynamics local to the grazing-sliding bifurcation
despite the fact that iterates of $\cP$ with $X > 0$ are virtual.

Next we compute $\cP$ to leading order.
The calculations are relatively routine and so details are omitted for brevity.
Let
\begin{equation}
A = \begin{bmatrix}
0 & 1 & 0 \\
0 & 0 & 1 \\
-\alpha_1 & -\alpha_2 & -\alpha_3
\end{bmatrix},
\label{eq:A}
\end{equation}
denote the Jacobian of the left half-system of \eqref{eq:odeEx}.
Let $\varphi_t(\bX_0,t_0)$ denote the solution to left half-system
with the arbitrary initial condition $\bX = \bX_0$ at $t = t_0$.
We have
\begin{equation}
\varphi_t(\bX_0,t_0) = \bX_p(t) + \bX_h(t;\bX_0,t_0),
\label{eq:generalSoln}
\end{equation}
where $\bX_p(t)$ is the particular solution \eqref{eq:Xp} and
\begin{equation}
\bX_h(t;\bX_0,t_0) = \re^{(t-t_0) A} \left( \bX_0 - \bX_p(t_0) \right),
\label{eq:Xh}
\end{equation}
is the homogeneous solution.
Via straight-forward but lengthy calculations using \eqref{eq:generalSoln}, we obtain
\begin{align}
\cP_g(X,t,Z) &=
\begin{bmatrix} 0 \\ t_{\rm graz} + 2 \pi \\ Z_{\rm graz} \end{bmatrix} +
\re^{2 \pi A} \begin{bmatrix} X \\ t - t_{\rm graz} \\ Z - Z_{\rm graz} \end{bmatrix} +
\frac{1}{\gamma_{\rm graz}} \left( I - \re^{2 \pi A} \right)
\begin{bmatrix} 1 \\ 0 \\ -1 \end{bmatrix}
\left( \gamma - \gamma_{\rm graz} \right) \nonumber \\
&\quad+
\cO \left( \left( X, t-t_{\rm graz}, Z-Z_{\rm graz}, \gamma-\gamma_{\rm graz} \right)^2 \right),
\label{eq:Pg}
\end{align}
where $Z_{\rm graz} = -1$, see \eqref{eq:Xgraz}.
The matrix part of $\cP_g$ has the particularly simple form $\re^{2 \pi A}$
due in part our choice of the ordering $(X,t,Z)$.
By using \eqref{eq:generalSoln} and the equations governing sliding motion \eqref{eq:slidingODE}, we also obtain
\begin{align}
\cP_d(X,t,Z) &=
\begin{bmatrix} 0 \\ t_{\rm graz} \\ Z_{\rm graz} \end{bmatrix} +
\begin{bmatrix}
0 & 0 & 0 \\
\beta_1 + 1 & 1 & 0 \\
\beta_2 & 0 & 1
\end{bmatrix}
\begin{bmatrix} X \\ t - t_{\rm graz} \\ Z - Z_{\rm graz} \end{bmatrix} \nonumber \\
&\quad+
X \cO \left( \sqrt{X}, t-t_{\rm graz}, Z-Z_{\rm graz}, \gamma-\gamma_{\rm graz} \right),
\label{eq:Pd}
\end{align}
for $X > 0$.
Refer to \cite{DiBu08,DiKo02,GlKo12} for detailed calculations of such a discontinuity map.

By then writing $\bx = (X,(t-t_{\rm graz}) {\rm ~mod~} 2 \pi,Z-Z_{\rm graz})$ and $\mu = \gamma-\gamma_{\rm graz}$,
to leading order $\cP$ is given by \eqref{eq:f} with
\begin{align}
A_L &= \re^{2 \pi A} \;, \label{eq:AL} \\
A_R &= \re^{2 \pi A}
\begin{bmatrix}
0 & 0 & 0 \\
\beta_1 + 1 & 1 & 0 \\
\beta_2 & 0 & 0
\end{bmatrix}, \label{eq:AR} \\
b &= \frac{1}{\gamma_{\rm graz}} \left( I - \re^{2 \pi A} \right)
\begin{bmatrix} 1 \\ 0 \\ -1 \end{bmatrix}. \label{eq:b}
\end{align}

\section{Fitting the parameters of the ODE system}
\label{sec:parameters}
\setcounter{equation}{0}

Here we determine values of
$\alpha_1$, $\alpha_2$, $\alpha_3$, $\beta_1$, and $\beta_2$
for which $A_L$ and $A_R$, as given by \eqref{eq:AL} and \eqref{eq:AR},
have desired sets of eigenvalues.

Let $\lambda^J_j$ denote the eigenvalues of $A_J$, for $j = 1,2,3$ and $J = L,R$,
with $\lambda^R_3 = 0$.
Since $A$, given by \eqref{eq:A}, is a real-valued matrix,
the eigenvalues of $A_L = \re^{2 \pi A}$ are either real and positive or appear in complex conjugate pairs.
Here we suppose that $\lambda^L_1 > 0$ and
$\lambda^L_{2,3} = p \pm \ri q$, for some $p \in \mathbb{R}$ and $q > 0$.
Then the eigenvalues of $A$ are
\begin{equation}
\begin{split}
\nu_1 &= \frac{1}{2 \pi} \ln \left( \lambda^L_1 \right), \\
\nu_{2,3} &= \frac{1}{4 \pi} \ln \left( p^2 + q^2 \right) \pm
\frac{\ri}{2 \pi} \tan^{-1} \left( \frac{q}{p} \right).
\end{split}
\label{eq:nu123}
\end{equation}
The trace, second trace, and determinant of $A$
are given by $-\alpha_3$, $\alpha_2$, and $-\alpha_1$ respectively, thus
the required values of $\alpha_1$, $\alpha_2$, and $\alpha_3$ are given by
\begin{equation}
\begin{split}
\alpha_1 &= -\nu_1 \nu_2 \nu_3 \;, \\
\alpha_2 &= \nu_1 \nu_2 + \nu_1 \nu_3 + \nu_2 \nu_3 \;, \\
\alpha_3 &= -(\nu_1 + \nu_2 + \nu_3),
\end{split}
\label{eq:gammai}
\end{equation}
see Appendix \ref{app:secondTrace}.

It remains for us to determine $\beta_1$ and $\beta_2$ in terms of $\lambda^R_1$ and $\lambda^R_2$.
Let $a_{ij}$ denote the $(i,j)$-element of $\re^{2 \pi A}$, for $i,j = 1,2,3$.
By using \eqref{eq:AR} to evaluate $\det \left( \lambda I - A_R \right)$, we obtain
\begin{align*}
0 &= \lambda^{R^2}_j - \left( a_{12} (\beta_1 + 1) + a_{13} \beta_2 + a_{22} + a_{33} \right) \lambda^R_j +
\left( a_{12} a_{33} - a_{13} a_{32} \right) (\beta_1 + 1) \nonumber \\
&\quad+
\left( a_{13} a_{22} - a_{12} a_{23} \right) \beta_2 +
a_{22} a_{33} - a_{23} a_{32} \;,
\end{align*}
for $j = 1,2$.
This provides two linear equations for $\beta_1$ and $\beta_2$, the solution to which is
\begin{equation}
\begin{split}
\beta_1 &= -1 + \frac{a_{12} a_{23} \left( \lambda^R_1 + \lambda^R_2 - a_{22} - a_{33} \right) +
a_{13} \left( \lambda^R_1 \lambda^R_2 - a_{22} \left( \lambda^R_1 + \lambda^R_2 \right) +
a_{23} a_{32} + a_{22}^2 \right)}
{a_{12}^2 a_{23} - a_{13}^2 a_{32} +
a_{12} a_{13} \left( a_{33} - a_{22} \right)} \;, \\
\beta_2 &= -\frac{a_{12} \left( \lambda^R_1 \lambda^R_2 - a_{33} \left( \lambda^R_1 + \lambda^R_2 \right) +
a_{23} a_{32} + a_{33}^2 \right) +
a_{13} a_{32} \left(\lambda^R_1 + \lambda^R_2 - a_{22} - a_{33} \right)}
{a_{12}^2 a_{23} - a_{13}^2 a_{32} +
a_{12} a_{13} \left( a_{33} - a_{22} \right)} \;,
\end{split}
\label{eq:betai}
\end{equation}
assuming $a_{12}^2 a_{23} - a_{13}^2 a_{32} + a_{12} a_{13} \left( a_{33} - a_{22} \right) \ne 0$,
as is generically the case\removableFootnote{
I do not see a simple geometric or physical interpretation for the case that this quantity is zero.
}.

\section{A bifurcation diagram}
\label{sec:bifDiag}
\setcounter{equation}{0}

Here we apply the formulas of \S\ref{sec:parameters} to the example of \S\ref{sub:verificationMainExample}.
This example is for the border-collision normal form with $\mu = 1$.
The eigenvalues of $A_L$, given by \eqref{eq:bcNormalForm},
are of the form $\lambda^L_1 > 0$ and $\lambda^L_{2,3} = p \pm \ri q$
for all points $(\sigma_L,\sigma_R) \in \cD$ that lie to the left of the dashed curve shown in Fig.~\ref{fig:domainExF}
and with $\sigma_R < 2.97$ approximately\removableFootnote{
$\sigma_R \approx 2.9656$ on the dashed curve.
}.

With the specific values \eqref{eq:exF}, corresponding to the black dot in Fig.~\ref{fig:domainExF},
the eigenvalues of $A_L$ and $A_R$ are\removableFootnote{
See {\sc goEx6all.m} and {\sc construct4dFilippov.m}.
}
\begin{align*}
\lambda^L_1 &\approx 0.2262333771 \;, \\
\lambda^L_{2,3} &\approx -0.3445852200 \pm 0.4870055259\ri \;, \\
\lambda^R_1 &= -1 \;, \\
\lambda^R_2 &= -1.75 \;, \\
\lambda^R_3 &= 0 \;,
\end{align*}
where each $\lambda^L_j$ is given to ten decimal places.
By substituting these values into \eqref{eq:nu123}--\eqref{eq:betai}, we obtain
\begin{equation}
\begin{split}
\alpha_1 &\approx 0.0302445699 \;, \\
\alpha_2 &\approx 0.1667559781 \;, \\
\alpha_3 &\approx 0.4009520660 \;, \\
\beta_1 &\approx -0.3783802961 \;, \\
\beta_2 &\approx -0.5981255840 \;,
\end{split}
\label{eq:param2b}
\end{equation}
to ten decimal places.

Before we discuss the dynamics of the ODE system \eqref{eq:odeEx} with the values \eqref{eq:param2b},
we first note that for the grazing-sliding bifurcation at $\gamma = \gamma_{\rm graz} \approx 0.9120$,
the return map $\cP$ is given by \eqref{eq:f} with \eqref{eq:AL}--\eqref{eq:b}, to leading order.
For this map we have $\det(\cO_L) \approx -5.4366$ and $\varrho^{\sf T} b \approx 1.7351$.
As discussed at the beginning of \S\ref{sec:derivingExamples},
since these quantities are nonzero the return map is conjugate to the 
border-collision normal form for $\mu \ne 0$.

With the given parameter values,
the border-collision normal form has infinitely many admissible, asymptotically stable $\cX^k \cY$-cycles
for $\mu > 0$ (Theorem \ref{th:twoParamEx}).
By conjugacy, the leading order approximation to $\cP$
has infinitely many admissible, asymptotically stable $\cX^k \cY$-cycles for $\gamma > \gamma_{\rm graz}$.
Furthermore, for each $k \ge k_{\rm min}$, the $\cX^k \cY$-cycle is a structurally stable invariant set.
Hence there exists $\gamma_k > \gamma_{\rm graz}$ such that $\cP$ 
has an admissible, asymptotically stable $\cX^k \cY$-cycle for all $\gamma \in (\gamma_{\rm graz},\gamma_k)$.

\begin{figure}[b!]
\begin{center}
\setlength{\unitlength}{1cm}
\begin{picture}(16.25,6.5)
\put(0,0){\includegraphics[height=6.5cm]{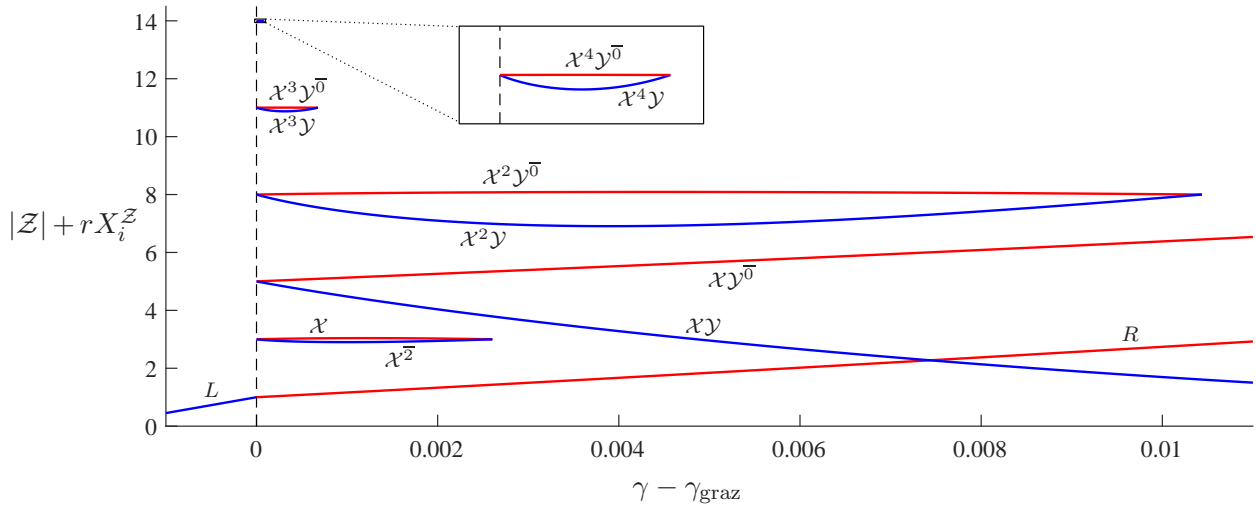}}
\put(8,0){\small $\gamma - \gamma_{\rm graz}$}
\put(-.3,3.48){\small $|\cZ| \hspace{-.2mm}+\hspace{-.2mm} r X^\cZ_i$}
\put(2.3,1.3){\scriptsize $L$}
\put(14.5,2.04){\scriptsize $R$}
\put(4.7,1.7){\scriptsize $\cX^{\overline{2}}$}
\put(3.7,2.16){\scriptsize $\cX$}
\put(8.7,2.16){\scriptsize $\cX \cY$}
\put(9,2.77){\scriptsize $\cX \cY^{\overline{0}}$}
\put(5.7,3.31){\scriptsize $\cX^2 \cY$}
\put(6,4.14){\scriptsize $\cX^2 \cY^{\overline{0}}$}
\put(3.16,4.81){\scriptsize $\cX^3 \cY$}
\put(3.16,5.26){\scriptsize $\cX^3 \cY^{\overline{0}}$}
\put(7.8,5.17){\scriptsize $\cX^4 \cY$}
\put(7.1,5.7){\scriptsize $\cX^4 \cY^{\overline{0}}$}
\end{picture}
\caption{
A bifurcation diagram of the ODE system \eqref{eq:odeEx} with \eqref{eq:param2b}.
Infinitely many asymptotically stable periodic orbits are created at $\gamma = \gamma_{\rm graz}$.
These correspond to the $\cX^k \cY$-cycles of Fig.~\ref{fig:qqExdRzero6}
and are shown here for $k = 1,\ldots,4$.
Curves corresponding to stable [unstable] periodic orbits are coloured blue [red].
\label{fig:bifDiagGrazSliding}
}
\end{center}
\end{figure}

We conclude that for the ODE system \eqref{eq:odeEx} with \eqref{eq:param2b},
infinitely many asymptotically stable periodic orbits
are created in the grazing-sliding bifurcation at $\gamma_{\rm graz}$.
For this example, $\cX$ and $\cY$ are words of length three and two, respectively.
Thus each $\cX^k \cY$-cycle of $\cP$ corresponds to a periodic orbit of \eqref{eq:odeEx}
that consists of $3 k + 2$ loops near the base periodic orbit $\bX_p(t)$, \eqref{eq:Xp}.
Since $\cX$ has two $R$'s and $\cY$ has one $R$,
exactly $2 k + 1$ of these loops involve a segment of sliding motion.

Fig.~\ref{fig:bifDiagGrazSliding} is a numerically computed bifurcation diagram of \eqref{eq:odeEx} with \eqref{eq:param2b}
illustrating several stable (blue) and unstable (red) periodic orbits.
Let us first explain the quantity $|\cZ| + r X^\cZ_i$ plotted on the vertical axis.
Each periodic orbit corresponds to a $\cZ$-cycle, for some $\cZ$, in the return map $\cP$ (e.g.~$\cZ = \cX^k \cY$).
We let $|\cZ|$ denote the length of $\cZ$
(this is also the number of loops that the periodic orbit has near the base periodic orbit).
For each $\cZ$ we choose a convenient index $i$ and let $X^\cZ_i = e_1^{\sf T} \bx^\cZ_i$ denote the first component of
the $i^{\rm th}$ point of the $\cZ$-cycle.
Finally, $r = 500$ is a scaling factor that enables the various periodic orbits to be distinguished clearly.
Note that each curve in Fig.~\ref{fig:bifDiagGrazSliding} has the integer value $|\cZ|$ at $\gamma = \gamma_{\rm graz}$
because as $\gamma \to \gamma_{\rm graz}^+$ the $\cZ$-cycle collapses to the origin and so here $X^\cZ_i = 0$.

Fig.~\ref{fig:bifDiagGrazSliding} shows periodic orbits corresponding
to $\cX^k \cY$-cycles and $\cX^k \cY^{\overline{0}}$-cycles for $k = 1,\ldots,4$.
For each $k$, these exist for all $\gamma \in (\gamma_{\rm graz},\gamma_k)$.
At $\gamma = \gamma_k$, the two periodic orbits collide and annihilate in a secondary
grazing-sliding bifurcation that mimics a saddle-node bifurcation.
In Fig.~\ref{fig:bifDiagGrazSliding} the two corresponding bifurcation curves intersect at $\gamma = \gamma_k$
because for each periodic orbit the index $i$ used for the vertical axis was chosen so that $X^\cZ_i = 0$ at $\gamma = \gamma_k$.
The values of $\gamma_k$ decrease as $k$ increases.
A determination of the asymptotic rate at which $\gamma_k \to 0$ as $k \to \infty$,
akin to that achieved in \cite{Si14b}, is beyond the scope of this paper.

Fig.~\ref{fig:bifDiagGrazSliding} shows the periodic orbit $\bX_p(t)$ (labelled $L$)
which exists for $\gamma < \gamma_{\rm graz}$.
For $\gamma > \gamma_{\rm graz}$ there exists one periodic orbit with one loop.
This periodic is unstable and involves a sliding segment (so is labelled $R$).
Also the periodic orbit corresponding to the $\cX$-cycle
exists for $\gamma \in (\gamma_{\rm graz},\gamma_{\rm graz} + 0.0026)$, approximately.
At the right end-point of this interval this periodic orbit
collides and annihilates with a periodic orbit corresponding to an $\cX^{\overline{2}} = RLL$-cycle.

The numerical continuation used to compute Fig.~\ref{fig:bifDiagGrazSliding}
was achieved by evaluating $\cP_g$ and $\cP_d$ numerically.
Computation of $\cP_g$ did not require a numerical ODE solver
because the exact solution is given by \eqref{eq:generalSoln}.
A numerical ODE solver was used to simulate sliding motion
governed by the nonlinear system \eqref{eq:slidingODE}.

Newton's method was used to locate fixed points of $\cP$.
This was achieved using the return map on $\Gamma$, call it $\tilde{\cP}$.
As described in \cite{GlJe15}, the map $\tilde{\cP}$
has the numerical advantage of being of one less dimension than that of $\cP$.
We did not discuss $\tilde{\cP}$ in \S\ref{sec:odeExample} because it has the analytical disadvantage
that each iterate of $\tilde{\cP}$ corresponds to several loops near $\bX_p(t)$
(specifically $\tilde{\cP} = \cP_d \circ \cP_g^\ell$,
where $\ell \ge 1$ is the number of loops required to reintersect $X = 0$).

Numerically continuing periodic orbits corresponding to $\cX^k \cY^{\overline{0}}$-cycles required particularly high precision,
not because they are unstable, but because they involve one point relatively close to the switching manifold.
On the switching manifold $\tilde{\cP}$ is non-differentiable,
so an extremely small discretisation was required to accurately estimate derivatives of $\tilde{\cP}$
for Newton's method.

\section{Discussion}
\label{sec:conc}
\setcounter{equation}{0}

The existence of multiple attractors in a dynamical system is a critical cause for complexity.
Here the long-term dynamics depends on the initial conditions
and in the presence of noise solutions may flip-flop between neighbourhoods of attractors.
A wide range of systems have been found to have large numbers of attractors
with varying physical consequences \cite{Fe08}.
For instance, a neuron that can exhibit a wide variety of stable bursting and beating solutions
appears to have the potential for sophisticated information processing \cite{CaBa93}.

At a grazing-sliding bifurcation,
an asymptotically stable periodic orbit can split into multiple attractors.
Since the attractors coincide at the bifurcation,
we cannot expect to know which attractor a particular orbit will converge to
if the parameter governing the bifurcation is varied dynamically \cite{DuNu99}.

This paper reveals that there is no limit to the number of attractors that can be created in a grazing-sliding bifurcation.
Infinitely many attractors are created for the example shown in Fig.~\ref{fig:bifDiagGrazSliding}.
These are destroyed in subsequent bifurcations shortly thereafter,
but there is no reason to expect that in other instances several attractors
cannot coexist over a relatively large region of parameter space.
The results have been demonstrated for an abstract ODE system
but are anticipated to occur in diverse physical systems due to the generality of the phenomenon.

\appendix

\section{Derivation of \eqref{eq:construct1}--\eqref{eq:construct3}}
\label{app:construction}
\setcounter{equation}{0}

Let $\lambda_1$, $\lambda_2$, and $\lambda_3$ be the eigenvalues of $M_\cX$.
The second trace of $M_\cX$ is $\sigma_\cX = \lambda_1 \lambda_2 + \lambda_1 \lambda_3 + \lambda_2 \lambda_3$,
see Appendix \ref{app:secondTrace}.
But $\cX$ includes at least one $R$,
thus the product $M_\cX = A_{\cX_{n-1}} \cdots A_{\cX_0}$ includes at least one instance of the matrix $A_R$.
Thus $M_\cX$ has a zero eigenvalue, say $\lambda_3$, because $\det(A_R) = \delta_R = 0$.
Thus to have $\lambda_2 = \frac{1}{\lambda_1}$ we require $\sigma_\cX = 1$.

If all conditions of Theorem \ref{th:SiTu17} are satisfied,
then the eigenvalues of $M_\cX$ are distinct (because $\lambda_1 > 1$, $\lambda_2 = \frac{1}{\lambda_1}$, and $\lambda_3 = 0$).
Thus there exist linearly independent $\zeta_1, \zeta_2, \zeta_3 \in \mathbb{R}^3$
with $M_\cX \zeta_i = \lambda_i \zeta_i$ for each $i$.
Since $\cX$ ends in $R$, we have
\begin{equation}
e_3^{\sf T} M_\cX \bx = 0 \;, \quad {\rm for~all~} \bx \in \mathbb{R}^3 \;.
\label{eq:MXx}
\end{equation}
Thus $e_3^{\sf T} \zeta_1 = 0$, because $e_3^{\sf T} \zeta_1 \lambda_1 = e_3^{\sf T} M_\cX \zeta_1 = 0$,
and similarly $e_3^{\sf T} \zeta_2 = 0$.
Then $e_3^{\sf T} \zeta_3 \ne 0$ by linear independence.

In view of the formula $\by_0 = \bx^\cX_0 - \frac{e_1^{\sf T} \bx^\cX_0}{e_1^{\sf T} \zeta_1} \,\zeta_1$ of Theorem \ref{th:SiTu17},
the vector $\bx^\cX_0 - \by_0$ is a scalar multiple of $\zeta_1$.
This implies that $\psi_1 = (I - M_\cX) \left( \bx^\cX_0 - \by_0 \right)$ is also a scalar multiple of $\zeta_1$,
where this expression is equivalent to \eqref{eq:psi1}
because $\bx^\cX_0$ is a fixed point of $f_\cX$ and so 
$\bx^\cX_0 = \left( I - M_\cX \right)^{-1} P_\cX b \mu$, by \eqref{eq:fX2}.
Therefore $M_\cX \psi_1 = \lambda_1 \psi_1$, and by using
\begin{equation}
\lambda_1 = \frac{e_1^{\sf T} M_\cX \psi_1}{e_1^{\sf T} \psi_1} \;,
\label{eq:lambda1}
\end{equation}
we see that $\xi_1$, given by \eqref{eq:xi1}, must be the zero vector.
The first component of $\xi_1$ is zero trivially;
the third component of $\xi_1$ is zero by \eqref{eq:MXx}.
For this reason we use the requirement that the second component of $\xi_1$ is zero
as our second condition on the parameter values.

As explained in \cite{SiTu17},
the conditions of Theorem \ref{th:SiTu17} imply that $M_\cY \zeta_1$
belongs to the stable subspace of $\bx^\cX_0$.
Since $\psi_1$ is a scalar multiple of $\zeta_1$,
the same is true for $\psi_2 = M_\cY \psi_1$.
Hence $\psi_2$ is a linear combination of $\zeta_2$ and $\zeta_3$.
But $e_3^{\sf T} \psi_2 = 0$, because $\cY$ ends in $R$.
Since $e_3^{\sf T} \zeta_2 = 0$ and $e_3^{\sf T} \zeta_3 \ne 0$,
the vector $\psi_2$ must be a scalar multiple of $\zeta_2$.
Thus $M_\cX \psi_2 = \frac{1}{\lambda_1} \psi_2$,
and by again using \eqref{eq:lambda1} we see that
$\xi_2$, given by \eqref{eq:xi2}, must be the zero vector.
We use the requirement that the first component of $\xi_2$ is zero
as our third condition on the parameter values.

\section{Proof of Theorem \ref{th:twoParamEx}}
\label{app:proof}
\setcounter{equation}{0}

Here we prove Theorem \ref{th:twoParamEx} by verifying
conditions (i)--(iv) of Theorem \ref{th:SiTu17}\removableFootnote{
Computations done in {\sc symExF2.m}.
I have tried to write my proof in such a way
that it does not appear to rely on computer arithmetic.
}:
\begin{enumerate}
\item
By directly substituting \eqref{eq:tRex1}--\eqref{eq:tLex1b} into $A_L$ and $A_R$,
we find that the matrix $M_\cX = A_R A_L A_R$ is given by
\begin{equation}
M_\cX = \begin{bmatrix}
\frac{(\sigma_R+1)^2}{\sigma_R^2+1} & \frac{\sigma_R^3-1}{\sigma_R^2+1} & -(\sigma_R+1) \\
\frac{\sigma_R^3-1}{\sigma_R \left( \sigma_R^2+1 \right)} &
\frac{\sigma_R^4 - \sigma_R^3 + \sigma_R^2 - \sigma_R + 1}{\sigma_R \left( \sigma_R^2+1 \right)} & -\sigma_R \\
0 & 0 & 0 
\end{bmatrix}.
\label{eq:MXex}
\end{equation}
The eigenvalues of this matrix are
\begin{equation}
\lambda_1 = \sigma_R + \frac{1}{\sigma_R} \;, \qquad
\lambda_2 = \frac{\sigma_R}{\sigma_R^2+1} \;, \qquad
\lambda_3 = 0 \;.
\label{eq:lambda123}
\end{equation}
Observe that $\lambda_2 = \frac{1}{\lambda_1}$ and $\lambda_1 > 1$ because $\sigma_R > 1$,
thus part (i) of Theorem \ref{th:SiTu17} is satisfied.
\item
The vectors
\begin{equation}
\zeta_1 = \begin{bmatrix} \sigma_R \\ \sigma_R-1 \\ 0 \end{bmatrix}, \qquad
\zeta_2 = \begin{bmatrix} -(\sigma_R-1) \\ 1 \\ 0 \end{bmatrix},
\label{eq:zeta12}
\end{equation}
are right eigenvectors of $M_\cX$ corresponding to $\lambda_1$ and $\lambda_2$.
The vectors 
\begin{align*}
\omega_1 &= \frac{1}{\sigma_R^2 - \sigma_R + 1} \begin{bmatrix} 1 \\ \sigma_R-1 \\ -\sigma_R \end{bmatrix}, \\
\omega_2 &= \frac{1}{\sigma_R^2 - \sigma_R + 1} \begin{bmatrix} -(\sigma_R-1) \\
\sigma_R \\ -\left( \sigma_R + \frac{1}{\sigma_R} \right) \end{bmatrix},
\end{align*}
are the corresponding left eigenvectors of $M_\cX$
normalised by $\omega_j^{\sf T} \zeta_j = 1$, for $j = 1,2$.
By directly evaluating $C = \begin{bmatrix} \omega_1^{\sf T} \\ \omega_2^{\sf T} \end{bmatrix} M_\cY
\begin{bmatrix} \zeta_1 & \zeta_2 \end{bmatrix}$, where $M_\cY = A_R A_L$, we obtain
(after much simplification)
\begin{equation}
\det(C) = \frac{1}{\sigma_R} \;.
\label{eq:c2}
\end{equation}
Part (ii) of Theorem \ref{th:SiTu17} is therefore satisfied
because $e_1^{\sf T} \zeta_1 \ne 0$ and $\lambda_2 < \det(C) < 1$
(due in part to $\sigma_R > 1$).
\item
The point $\bx^\cX_0$ is a fixed point of $f_\cX$,
thus $\bx^\cX_0 = \left( I - M_\cX \right)^{-1} P_\cX b \mu$, by \eqref{eq:fX2}.
With $\cX = RLR$ and $b \mu = e_1$, we have
$\bx^\cX_0 = \left( I - A_R A_L A_R \right)^{-1} \left( I + A_R + A_R A_L \right) e_1$.
By evaluating this expression and using
$\bx^\cX_1 = f_R \left( \bx^\cX_0 \right)$ and 
$\bx^\cX_2 = f_L \left( \bx^\cX_1 \right)$, we obtain
\begin{align}
\bx^\cX_0 &= \frac{1}{\sigma_R^2 - \sigma_R + 1}
\begin{bmatrix} \sigma_R^2 \\ -1 \\ 0 \end{bmatrix}, \\
\bx^\cX_1 &= \frac{1}{\sigma_R^2 - \sigma_R + 1}
\begin{bmatrix} -\sigma_R \left( \sigma_R^2+1 \right) \\ -\sigma_R^3 \\ 0 \end{bmatrix}, \\
\bx^\cX_2 &= \frac{1}{\sigma_R^3 + 1}
\begin{bmatrix} \sigma_L \sigma_R \left( \sigma_R^2+1 \right) - (\sigma_R-1) \\
\sigma_L \sigma_R (\sigma_R + 1) \left( \sigma_R^2 + 1 \right) \\
\left( \sigma_L \sigma_R^2 - 1 \right) \left( \sigma_R^2+1 \right) \end{bmatrix}.
\label{eq:xX2}
\end{align}
Since $\sigma_R > 1$, we have $e_1^{\sf T} \bx^\cX_0 > 0$ and $e_1^{\sf T} \bx^\cX_1 < 0$.
Also $e_1^{\sf T} \bx^\cX_2 > 0$ because
$\sigma_L > \frac{\sigma_R - 1}{\sigma_R \left( \sigma_R^2 + 1 \right)}$.
Thus the $\cX$-cycle is admissible with no points on $\Sigma$.
\item
As noted in \S\ref{sub:derivationMainExample}, we have $\by_0 = [0,-1,0]^{\sf T}$.
From the above expressions it is evident that this agrees with the formula
$\by_0 = \bx^\cX_0 - \frac{e_1^{\sf T} \bx^\cX_0}{e_1^{\sf T} \zeta_1} \,\zeta_1$.
By iterating $\by_0$ under $f$ we obtain $\by_1 = [0,0,0]^{\sf T}$ and $\by_2 = [1,0,0]^{\sf T}$.
Notice $\by_\alpha \in \Sigma$, as $\alpha = 1$.

Since $\by_0 - \bx^\cX_0$ is a scalar multiple of $\zeta_1$,
the point $\by_0$ has a backwards orbit that converges to the $\cX$-cycle
on the open line segments connecting $\by_i$ and $\bx^\cX_i$, for $i = 0,1,2$.
Moreover, $s(\by_i) = \cS_i$ for all $i \le 0$
because these line segments do not intersect $\Sigma$.

Similarly, from \eqref{eq:xX2} we find that $\by_2 - \bx^\cX_0$ is a scalar multiple of $\zeta_2$.
Thus the forward orbit of $\by_2$ converges to the $\cX$-cycle
on the open line segments connecting $\by_{(i+2) \,{\rm mod}\, n}$ and $\bx^\cX_i$, for $i = 0,1,2$,
and so $s(\by_i) = \cS_i$ for all $i \ge 2$,
This establishes that part (iv) of Theorem \ref{th:SiTu17} is satisfied.
\end{enumerate}
Thus all conditions of Theorem \ref{th:SiTu17} are satisfied,
hence there exists $k_{\rm min} \ge 0$ such that \eqref{eq:f} has an
admissible, asymptotically stable $\cX^k \cY$-cycle with no points on $\Sigma$ for all $k \ge k_{\rm min}$.
\hfill $\Box$

\section{The characteristic polynomial of a $3 \times 3$ matrix}
\label{app:secondTrace}
\setcounter{equation}{0}

Throughout this paper the second trace of a $3 \times 3$ matrix is particularly important.
Since the second trace arises less commonly than the ubiquitous trace and determinant,
here we state the basic properties of the characteristic polynomial of a $3 \times 3$ matrix
which includes the second trace.

The characteristic polynomial of a $3 \times 3$ matrix $Q$ is
\begin{equation}
\det(\lambda I - Q) = \lambda^3 - \tau \lambda^2 + \sigma \lambda - \delta \;,
\label{eq:charPoly}
\end{equation}
where $\tau$, $\sigma$, and $\delta$ 
are the trace, second trace, and determinant of $Q$, respectively.
Let $q_{ij}$ denote the $(i,j)$-element of $Q$, for $i,j = 1,2,3$.
By directly evaluating \eqref{eq:charPoly} we obtain
\begin{align*}
\tau &= q_{11} + q_{22} + q_{33} \;, \\
\sigma &= q_{11} q_{22} + q_{11} q_{33} + q_{22} q_{33} - q_{12} q_{21} - q_{13} q_{31} - q_{23} q_{32} \;, \\
\delta &= q_{11} q_{22} q_{33} + q_{12} q_{23} q_{31} + q_{13} q_{21} q_{32}
- q_{11} q_{23} q_{32} - q_{12} q_{21} q_{33} - q_{13} q_{22} q_{31} \;.
\end{align*}
Let $\lambda_1$, $\lambda_2$, and $\lambda_3$ be the eigenvalues of $Q$ (counting multiplicity).
Then $\det(\lambda I - Q) = (\lambda - \lambda_1)(\lambda - \lambda_2)(\lambda - \lambda_3)$,
and by matching this to \eqref{eq:charPoly} we obtain
\begin{align*}
\tau &= \lambda_1 + \lambda_2 + \lambda_3 \;, \\
\sigma &= \lambda_1 \lambda_2 + \lambda_1 \lambda_3 + \lambda_2 \lambda_3 \;, \\
\delta &= \lambda_1 \lambda_2 \lambda_3 \;.
\end{align*}

\section*{Acknowledgements}
The author thanks Mike Jeffrey for many useful discussions regarding this work.


\end{document}